\documentclass{amsproc}

\newtheorem{theorem}{Theorem}[section]

\newtheorem{proposition}[theorem]{Proposition}

\theoremstyle{definition}

\theoremstyle{remark}

\numberwithin{equation}{section}

\usepackage{mathrsfs}
\usepackage{enumitem}
\usepackage{xypic}

\newcommand{\vs}[0]{\vspace{2mm}}
\newcommand{\ol}[1]{\overline{#1}}
\newcommand{\ul}[1]{\underline{#1}}

\newcommand{\h}[0]{\varphi}

\usepackage{amssymb}

\newcommand{\smallmattwo}[4]{
\left(
\begin{smallmatrix}
#1 & #2 \\
#3 & #4 \\
\end{smallmatrix}
\right)
}

\begin{document}

\title[A trilogy of MCG representations from 3d quantum gravity]{A trilogy of mapping class group representations from three-dimensional quantum gravity}

\author{Hyun Kyu Kim}
\address{School of Mathematics, Korea Institute for Advanced Study (KIAS), 85 Hoegiro Dongdaemun-gu, Seoul 02455, Republic of Korea}
%\curraddr{}
\email{hkim@kias.re.kr}
%\thanks{}

%    General info
\date{March 23, 2023}

\dedicatory{To Igor B. Frenkel on the occasion of his 70th birthday}


\begin{abstract}
For a punctured surface $\mathfrak{S}$, the author and Scarinci (\cite{KS}) have recently constructed a quantization of a moduli space of Lorentzian metrics on the 3-manifold $\mathfrak{S} \times \mathbb{R}$ of constant sectional curvature $\Lambda \in \{-1,0,1\}$. The invariance of this quantization under the action of the mapping class group ${\rm MCG}(\mathfrak{S})$ of $\mathfrak{S}$ yields families of unitary representations of ${\rm MCG}(\mathfrak{S})$ on a Hilbert space, with key ingredients being three versions of the quantum dilogarithm functions depending on $\Lambda$. In this survey article, we review and elaborate on this result.
\end{abstract}

\maketitle

\setcounter{secnumdepth}{2}
\setcounter{tocdepth}{2}

\vspace{-5mm}

\tableofcontents

\section{Introduction}
\label{sec:introduction}

Quantization of Teichm\"uller spaces of surfaces appeared in 1990's as an approach to quantization of 3-dimensional gravity, as it had been expected in the mathematical physics literature since 1980's that the phase space of 3-dimensional pure gravity theory would be closely related to the Teichm\"uller space. Kashaev \cite{Kash98}, and independently Chekhov and Fock \cite{CF, F}, established the first major constructions of the quantum Teichm\"uller theory in late 1990's. What are quantized in these works are certain versions of the {\it Teichm\"uller space} of an oriented punctured surface $\frak{S}$. The Teichm\"uller space of $\frak{S}$, which dates back at least to early 20${}^{\rm th}$ century, can be viewed as the space of all complex structures on $\frak{S}$ considered up to pullback by self-diffeomorphisms of $\frak{S}$ isotopic to identity, and can also be viewed as the space of all complete hyperbolic metrics on $\frak{S}$ considered up to isotopy. There exist some variants considered in 1980's, depending on what kind of conditions or extra data one wants to put at the punctures. The Teichm\"uller space is equipped with a canonical Poisson structure called the Weil-Petersson Poisson structure, and there is a natural pullback action on the Teichm\"uller space, preserving the Weil-Petersson Poisson structure, by the {\it mapping class group} ${\rm MCG}(\frak{S})$ of $\frak{S}$, which is defined as the group of isotopy classes of orientation-preserving self-diffeomorphisms of $\frak{S}$. The above mentioned works establish certain versions of quantization of the Teichm\"uller space with respect to the Weil-Petersson Poisson structure, equivariant under the action of ${\rm MCG}(\frak{S})$. They build a complex Hilbert space $\mathscr{H}$, a set of self-adjoint operators on $\mathscr{H}$ that play roles of quantum counterparts of classical observable functions on the Teichm\"uller space, and a family of unitary projective representations of ${\rm MCG}(\frak{S})$ on $\mathscr{H}$ such that the set of (quantum) self-adjoint operators on $\mathscr{H}$ satisfy certain equivariance under these representations of ${\rm MCG}(\frak{S})$. It is often said that this family of representations of ${\rm MCG}(\frak{S})$ is a main consequence of the quantum Teichm\"uller theory. These representations are interpreted as forming a family because they depend on the real quantum parameter $\hbar$, usually called the Planck constant. 

\vs

These early works are based on special coordinate systems on the Teichm\"uller space developed in 1980's by Penner and Thurston \cite{P87, Thurston}, as well as on a special function called the {\it quantum dilogarithm} studied in 1990's by Faddeev and Kashaev \cite{Faddeev, FK94}. Especially, the framekwork of Chekhov and Fock \cite{F, CF} later generalized by Fock and Goncharov to the theory of quantization of cluster varieties \cite{FG09}, a prominent application being the quantization of higher Teichm\"uller spaces \cite{FG06}. So, the Chekhov-Fock quantization of the Teichm\"uller space can be viewed as a prototypical example of the Fock-Goncharov quantization of cluster varieties.

\vs

In the meantime, as the initial ideas of quantum Teichm\"uller theory were evolving to more general setting as mentioned above, the motivating original question of three dimensional quantum gravity seems to have been lost a bit in this process of generalization. Only relatively recently, Meusburger and Scarinci \cite{MS} found a direct interpretation of the moduli space of a certain version of the three dimensional gravity in terms of the language of cluster varieties. Building on this work, Scarinci and the author \cite{KS} developed a corresponding quantum theory, by suitably modifying Fock and Goncharov's quantum theory of cluster varieties. Roughly speaking, what Fock and Goncharov quantized in \cite{FG09} is the Poisson manifold given as the set of positive real points of a cluster variety \cite{FG06, FG09}. The work \cite{MS} realized a certain moduli space of 3d gravity as the set of points of a cluster variety valued at a special ring
\begin{align}
\label{eq:R_Lambda}
\mathbb{R}_\Lambda = \mathbb{R}[\ell]/(\ell^2=-\Lambda)
\end{align}
called the ring of {\it generalized complex numbers}, where $\Lambda$ is a real parameter that appears in the theory of gravity, called the {\it cosmological constant}. In \cite{KS} we developed an $\mathbb{R}_\Lambda$-version of the Fock-Goncharov quantization, and consequently were able to establish a version of the 3d quantum gravity. Just as quantum Teichm\"uller theory, one of the output of our 3d quantum gravity result \cite{KS} would be a family of unitary projective representations of ${\rm MCG}(\frak{S})$ on a Hilbert space, which this time depends not only on the quantum parameter $\hbar$, but also on the gravity parameter $\Lambda$. In fact, the formulation of \cite{KS} can be viewed as a construction of a family of representations of a certain groupoid version of ${\rm MCG}(\frak{S})$, but not exactly ${\rm MCG}(\frak{S})$. Hence, some extra work is needed if one wants to extract concrete ${\rm MCG}(\frak{S})$ representations from the results of \cite{KS}. In \cite{KS}, this task is not explicitly done, but only hints are given. 

\vs

We undertake this task in the present article. Namely, we review part of the 3d quantum gravity result of \cite{KS} together with some background materials, and provide necessary extra algebraic and functional-analytic treatments to construct unitary projective representations of ${\rm MCG}(\frak{S})$. We note that the way how the results of \cite{KS} are reviewed is adapted to the goal of the present article, and hence may contain somewhat different underlying viewpoints from those of \cite{KS}, which could provide new insights to the subject. We also touch upon some applications and possible future research topics.

\section{A three-dimensional quantum gravity}
\label{sec:overview}

\subsection{Basic definitions}

Throughout the present article, we let $\mathfrak{S}$ be an oriented punctured surface of finite type; one can view it as being obtained from the closed oriented surface of genus $g$ by removing $n$ points, called {\it punctures}, such that $n>1$ and $\chi(\mathfrak{S}) = 2 - 2g - n <0$. The moduli space of 3d gravity to be quantized is $\mathcal{GH}_\Lambda(\mathfrak{S} \times \mathbb{R})$ (see \cite{MS}), depending on the parameter $\Lambda \in \mathbb{R}$ called the cosmological constant. As only the sign of $\Lambda$ matters, it suffices to consider only three possible values $-1$, $0$ and $1$ for $\Lambda$. The space $\mathcal{GH}_\Lambda(\mathfrak{S} \times \mathbb{R})$ parametrizes isotopy classes of maximal globally hyperbolic Cauchy-complete Lorentzian metrics $g$ on the 3-manifold $\mathfrak{S} \times \mathbb{R}$ of constant sectional curvature $\Lambda$ having parabolic holonomy around punctures. 

\vs

Let us briefly explain the notions just appeared for the metrics $g$. A metric $g$ on $\mathfrak{S}\times \mathbb{R}$ of Lorentzian signature $(++-)$ being globally hyperbolic means that there is a Cauchy surface $\mathfrak{S}$, i.e. an embedded surface in $\mathfrak{S}\times \mathbb{R}$ diffeomorphic to $\mathfrak{S}$ that intersects each inextendible timelike curve exactly once. Such $g$ being maximal means that the Lorentzian manifold $(\mathfrak{S},g)$ is not properly contained via isometric embedding into a larger globally hyperbolic Lorentzian $3$-manifold. Such $g$ is Cauchy-complete if the induced metric on $\mathfrak{S}$ is geodesically complete. Such a Lorentzian manifold $(\mathfrak{S}\times \mathbb{R},g)$ can be obtained as a quotient of the model geometry, which is the anti-de Sitter space ${\rm AdS}^3$ for $\Lambda=-1$, the Minkowski space ${\rm Mink}^3$ for $\Lambda=0$, and the de Sitter space ${\rm dS}^3$ for $\Lambda=1$, by a discrete subgroup of the isometry group of the model geometry, which is ${\rm PSL}_2(\mathbb{R}_\Lambda)$, where $\mathbb{R}_\Lambda$ denotes the ring of generalized complex numbers $
\mathbb{R}_\Lambda = \mathbb{R}[\ell]/(\ell^2=-\Lambda)$ (see eq.\eqref{eq:R_Lambda}). Thus a point of $\mathcal{GH}_\Lambda(\frak{S}\times \mathbb{R})$ yields an element of ${\rm Hom}(\pi_1(\mathfrak{S}),{\rm PSL}_2(\mathbb{R}_\Lambda))/{\rm PSL}_2(\mathbb{R}_\Lambda)$, i.e. a group homomorphism from $\pi_1(\mathfrak{S})$ to ${\rm PSL}_2(\mathbb{R}_\Lambda)$, called the holonomy representation, defined up to conjugation in ${\rm PSL}_2(\mathbb{R}_\Lambda)$. The parabolicity condition is understood in terms of the holonomy.

\vs

An {\it ideal triangulation} $\Delta$ of $\mathfrak{S}$ means a collection of unoriented simple paths in $\mathfrak{S}$ called {\em ideal arcs} running between punctures, dividing $\mathfrak{S}$ into regions called {\em ideal triangles} bounded by three ideal arcs. In the present article we do not allow an ideal triangle to be `self-folded'; this means that we require the three sides of an ideal triangle to be mutually distinct ideal arcs. An ideal arc is often considered up to isotopy, and an ideal triangulation up to simultaneous isotopy of its constituent ideal arcs. It is a simple exercise to show that $\Delta$ has $6g-6+3n$ constituent ideal arcs, and forms $4g-4+2n$ ideal triangles.

\vs

Per each choice of an ideal triangulation $\Delta$ of $\mathfrak{S}$, Meusburger and Scarinci \cite{MS} constructed a special coordinate system on the space $\mathcal{GH}_\Lambda(\mathfrak{S} \times \mathbb{R})$, with the coordinate functions $x_i$ and $y_i$ associated to ideal arcs $i$ of $\Delta$. In fact, these do not form a genuine coordinate system, but a constrained one. Namely, these coordinates satisfy the constraint equations
\begin{align}
\label{eq:3d_constraint_equations}
\sum_{i \in \Delta} v_{i,p} x_i = 0, \qquad
\sum_{i \in \Delta} v_{i,p} y_i = 0
\end{align}
for each puncture $p$ of $\mathfrak{S}$, where $v_{i,p} \in \{0,1,2\}$ denotes the valence of $\Delta$ at $p$; we note that these constraints correspond to the parabolicity condition in the definition of $\mathcal{GH}_\Lambda(\mathfrak{S}\times \mathbb{R})$. What is useful to know is that for each puncture $p$, the vector $(v_{i,p})_{i\in \Delta}$, viewed as a column vector, belongs to the kernel of the exchange matrix $\varepsilon$, i.e.
\begin{align}
\label{eq:v_in_the_kernel}
\sum_{j\in \Delta} \varepsilon_{ij} v_{j,p} = 0, \qquad \forall i \in \Delta, \quad \forall \mbox{puncture } p.
\end{align}
It is well known that these vectors actually form a basis of the entire kernel of $\varepsilon$ (see e.g. \cite[Lem.8]{BL}). The natural gravitational Poisson structure $\{\cdot,\cdot\}$ of $\mathcal{GH}_\Lambda(\mathfrak{S} \times \mathbb{R})$ takes values on these coordinate functions as
$$
\{x_i, x_j\}= 0, \qquad
\{y_i,y_j\} = 0, \qquad
\{x_i, y_i \} = \varepsilon_{ij}, \qquad \forall i,j \in \Delta.
$$

\vs

We notice that there is a natural action on $\mathcal{GH}_\Lambda(\mathfrak{S}\times \mathbb{R})$ of the {\em mapping class group} ${\rm MCG}(\mathfrak{S}) = {\rm Diff}_+(\mathfrak{S})/{\rm Diff}(\mathfrak{S})_0$ of $\mathfrak{S}$, which is the group of orientation-preserving self-diffeomorphisms of $\mathfrak{S}$ considered up to isotopy, preserving the Poisson structure. See \cite{MS,KS} for more details and references for the structure of the moduli space $\mathcal{GH}_\Lambda(\mathfrak{S}\times \mathbb{R})$.

\subsection{A brief overview of a quantization program}
\label{subsec:brief_overview}

We describe the formulation and results of the quantization program established in the work of the author and Scarinci \cite{KS}. For each ideal triangulation $\Delta$, using the above mentioned coordinate system on $\mathcal{GH}_\Lambda(\frak{S}\times \mathbb{R})$ a certain class of classical observables $\mathcal{A}_\Delta = \mathcal{A}_{\Delta;\Lambda}$ is suggested, which depends on the cosmological constant $\Lambda \in \{-1,0,1\}$.  The best way to view this class is to formulate it as an algebra of certain functions on the space $\mathcal{GH}_\Lambda(\frak{S}\times \mathbb{R})$ valued in the ring of generalized complex numbers
$
\mathbb{R}_\Lambda = \mathbb{R}[\ell]/(\ell^2=-\Lambda)
$, 
or maybe in a suitable complexification of $\mathbb{R}_\Lambda$. In the present article we do not delve into a detailed discussion about it, and one can just view $\mathcal{A}_\Delta$ as a certain subset of $C^\infty(\mathcal{GH}_\Lambda(\frak{S}\times \mathbb{R}))$, which is a collection of functions to be quantized. A drawback that we have to endure by giving up working with $\mathbb{R}_\Lambda$ is that the class $\mathcal{A}_\Delta \subset C^\infty(\mathcal{GH}_\Lambda(\frak{S}\times \mathbb{R}))$ may not form an algebra. 

\vs

Next, a quantum counterpart $\hat{\mathcal{A}}_\Delta^\hbar$ is introduced, which could be viewed as a one-real-parameter family of non-commutative algebras that deforms $\mathcal{A}_\Delta$ in the direction of the Poisson structure. A representation $\pi^\hbar_\Delta$ of $\hat{\mathcal{A}}^\hbar_\Delta$ is constructed on a Hilbert space $\mathscr{H}_\Delta$, in an appropriate sense. In particular, to each element $\hat{u} \in \hat{\mathcal{A}}^\hbar_\Delta$, is associated a densely-defined operator $\pi^\hbar_\Delta(\hat{u})$ on $\mathscr{H}_\Delta$. There are some subtle functional-analytic issues on the domains of definition, but we do not discuss them in this article.

\vs

The classes $\mathcal{A}_\Delta$ and $\mathcal{A}^\hbar_\Delta$ are constructed in such a way that for each pair of ideal triangulations $\Delta$ and $\Delta'$, there are natural bijections of sets (or of algebras)
$$
\mu_{\Delta,\Delta'} : \mathcal{A}_{\Delta'} \to \mathcal{A}_\Delta
\quad\mbox{and}\quad
\mu^\hbar_{\Delta,\Delta'} : \mathcal{A}^\hbar_{\Delta'} \to \mathcal{A}^\hbar_\Delta,
$$
that satisfy the consistency equations
$$
\mu_{\Delta,\Delta'} \circ \mu_{\Delta',\Delta''} = \mu_{\Delta,\Delta''} \quad\mbox{and}\quad
\mu^\hbar_{\Delta,\Delta'} \circ \mu^\hbar_{\Delta',\Delta''} = \mu^\hbar_{\Delta,\Delta''}
$$
for each triple of ideal triangulations $\Delta$, $\Delta'$ and $\Delta''$. Hence, the classes $\mathcal{A}_\Delta$ for different $\Delta$ can be consistently identified with each other, and likewise for $\hat{\mathcal{A}}^\hbar_\Delta$.

\vs

For each pair of ideal triangulations $\Delta$ and $\Delta'$, a unitary map
$$
{\bf K}^\hbar_{\Delta,\Delta'} : \mathscr{H}_{\Delta'} \to \mathscr{H}_\Delta
$$
between the Hilbert spaces associated to $\Delta$ and $\Delta'$ is constructed, so that it intertwines the representations $\pi^\hbar_{\Delta'}$ and $\pi^\hbar_\Delta$ of $\hat{\mathcal{A}}^\hbar_{\Delta'}$ and $\hat{\mathcal{A}}^\hbar_\Delta$, which are related via the map $\mu^\hbar_{\Delta,\Delta'}$. To be more precise, for each $\hat{u} \in \hat{\mathcal{A}}^\hbar_{\Delta'}$ we would like the diagram
$$
\xymatrix@C+10mm{
\mathscr{H}_{\Delta'} \ar[r]^-{\pi^\hbar_{\Delta'}(\hat{u})} \ar[d]_{{\bf K}^\hbar_{\Delta,\Delta'}} & \mathscr{H}_{\Delta'} \ar[d]^{{\bf K}^\hbar_{\Delta,\Delta'}}  \\
\mathscr{H}_\Delta \ar[r]^-{\pi^\hbar_\Delta(\mu^\hbar_{\Delta,\Delta'}(\hat{u}))} & \mathscr{H}_\Delta
}
$$
to commute, that is, we would like the intertwining equations
$$
{\bf K}^\hbar_{\Delta,\Delta'} \circ \pi^\hbar_{\Delta'}(\hat{u}) = \pi^\hbar_{\Delta,\Delta'}(\mu^\hbar_{\Delta,\Delta'}(\hat{u})) \circ {\bf K}^\hbar_{\Delta,\Delta'}
$$
to hold. Another important stipulation is that these {\it unitary intertwiners} ${\bf K}^\hbar_{\Delta,\Delta'}$ should satisfy the consistency equations
$$
{\bf K}^\hbar_{\Delta,\Delta'} \circ {\bf K}^\hbar_{\Delta',\Delta''} = {\bf K}^\hbar_{\Delta,\Delta''}
$$
for each triple of triangulations $\Delta$, $\Delta'$ and $\Delta''$, perhaps up to multiplicative constants. So, the unitary maps ${\bf K}^\hbar_{\Delta,\Delta'}$ let us identify the representations $\pi^\hbar_\Delta$ of $\hat{\mathcal{A}}^\hbar_\Delta$ on the Hilbert spaces $\mathscr{H}_\Delta$ for different $\Delta$, in a consistent manner. 

\vs

An algebraic deformation quantization map
$$
\hat{Q}^\hbar_\Delta : \mathcal{A}_\Delta \to \hat{\mathcal{A}}^\hbar_\Delta
$$
is constructed for each $\Delta$. The most important property of these maps $\hat{Q}^\hbar_\Delta$ is the compatibility with respect to the algebraic compatibility maps $\mu_{\Delta,\Delta'}$ and $\mu^\hbar_{\Delta,\Delta'}$. Namely, they make the diagram
$$
\xymatrix@R-3mm{
\mathcal{A}_{\Delta'} \ar[r]^-{\hat{Q}^\hbar_{\Delta'}} \ar[d]_{\mu_{\Delta,\Delta'}} & \hat{\mathcal{A}}^\hbar_{\Delta'} \ar[d]^{\mu^\hbar_{\Delta,\Delta'}} \\
\mathcal{A}_\Delta \ar[r]^-{\hat{Q}^\hbar_\Delta} & \hat{\mathcal{A}}^\hbar_\Delta \\
}
$$
commute, that is, the following equation holds for each pair of ideal triangulations $\Delta$ and $\Delta'$:
$$
\mu^\hbar_{\Delta,\Delta'} \circ \hat{Q}^\hbar_{\Delta'} = \hat{Q}^\hbar_\Delta \circ \mu_{\Delta,\Delta'}.
$$
This means that the deformation quantization maps $\hat{Q}^\hbar_\Delta$ for different $\Delta$ can be identified with each other in a consistent manner. The construction of these maps $\hat{Q}^\hbar_\Delta$ is built on a previous joint work of the author and Allegretti \cite{AK}, which in turn is based on the work of Bonahon and Wong \cite{BW}.

\vs

Let us summarize. For each ideal triangulation $\Delta$ a quantization is constructed in the form of the following composition
$$
\xymatrix@C+5mm{
\mathcal{A}_\Delta \ar[r]^-{\hat{Q}^\hbar_\Delta} & \hat{\mathcal{A}}^\hbar_\Delta \ar[r]^-{\pi^\hbar_\Delta} & \{\mbox{densely-defined operators on $\mathscr{H}_\Delta$}\}
}.
$$
For each pair of ideal triangulations $\Delta$ and $\Delta'$, these quantizations are compatible with each other in the sense that the following diagram commutes:
$$
\xymatrix@R+1mm@C+5mm{
\mathcal{A}_{\Delta'} \ar[r]^-{\hat{Q}^\hbar_{\Delta'}} \ar[d]_{\mu_{\Delta,\Delta'}} & \hat{\mathcal{A}}^\hbar_{\Delta'} \ar[r]^-{\pi^\hbar_{\Delta'}} \ar[d]^{\mu^\hbar_{\Delta,\Delta'}} & \{\mbox{densely-defined operators on $\mathscr{H}_{\Delta'}$}\} \ar[d]^{\mbox{conjugation by ${\bf K}^\hbar_{\Delta,\Delta'}$}} \\
\mathcal{A}_\Delta \ar[r]^-{\hat{Q}^\hbar_\Delta} & \hat{\mathcal{A}}^\hbar_\Delta \ar[r]^-{\pi^\hbar_\Delta} & \{\mbox{densely-defined operators on $\mathscr{H}_\Delta$}\} \\
}
$$
Hence, in the end, one can say that the quantization does not depend on the choice of $\Delta$, and thus write the result as
$$
\xymatrix@C+5mm{
\mathcal{A}_\frak{S} \ar[r]^-{\hat{Q}^\hbar_\frak{S}} & \hat{\mathcal{A}}^\hbar_\frak{S} \ar[r]^-{\pi^\hbar_\frak{S}} & \{\mbox{densely-defined operators on $\mathscr{H}_\frak{S}$}\}
},
$$
which depends only on the surface $\frak{S}$.

\subsection{On the mapping class group action}

An important stipulation we should impose on a quantization is an invariance or equivariance under the mapping class group ${\rm MCG}(\frak{S})$, which acts naturally on the moduli space $\mathcal{GH}_\Lambda(\frak{S} \times \mathbb{R})$ as Poisson automorphisms, hence on $C^\infty(\mathcal{GH}_\Lambda(\frak{S}\times \mathbb{R}))$. 

\vs

We first require that the class of classical observables $\mathcal{A}_\Delta$ should be invariant under the action of ${\rm MCG}(\frak{S})$. That is, we should have $\mathcal{A}_\Delta = \mathcal{A}_{g.\Delta}$ for $g\in {\rm MCG}(\frak{S})$, where $g.\Delta$ denotes the isotopy class of an ideal triangulation obtained by acting the mapping class $g$ to $\Delta$; the bijections between these classical classes associated to changes of ideal triangulations must also satisfy the invariance $\mu_{\Delta,\Delta'} = \mu_{g.\Delta,g.\Delta'}$. Upon quantization, invariance should hold for the classes of quantum observables $\mathcal{A}^\hbar_\Delta = \mathcal{A}^\hbar_{g.\Delta}$, the quantum bijections between them $\mu^\hbar_{\Delta,\Delta'} = \mu^\hbar_{g.\Delta,g.\Delta'}$, the quantum Hilbert spaces $\mathscr{H}_\Delta = \mathscr{H}_{g.\Delta}$, the representations $\pi^\hbar_\Delta = \pi^\hbar_{g.\Delta}$, the unitary intertwiners between them ${\bf K}^\hbar_{\Delta,\Delta'} = {\bf K}^\hbar_{g.\Delta,g.\Delta'}$, and the deformation quantization maps $\hat{Q}^\hbar_\Delta = \hat{Q}^\hbar_{g.\Delta}$.

\vs

In fact, one should carefully think of the meaning of the equalities appearing in the invariance condition, such as $\mathcal{A}_\Delta = \mathcal{A}_{g.\Delta}$ and $\mathscr{H}_\Delta=\mathscr{H}_{g.\Delta}$. For simplicity of discussion, at the moment, assume that $\mathcal{A}_\Delta$ is realized as a subalgebra of $C^\infty(\mathcal{GH}_\Lambda(\frak{S}\times \mathbb{R}))$ for each $\Delta$; this is indeed the case when $\Lambda=-1$. Thus both $\mathcal{A}_\Delta$ and $\mathcal{A}_{g.\Delta}$ are subalgebras of $C^\infty(\mathcal{GH}_\Lambda(\frak{S}\times \mathbb{R}))$, hence the equality $\mathcal{A}_\Delta = \mathcal{A}_{g.\Delta}$ seems to have an unambiguous meaning. However, for a better control, we should discuss what is a correct identification map between $\mathcal{A}_\Delta$ and $\mathcal{A}_{g.\Delta}$. An obvious option is the restriction of the identity map on $C^\infty(\mathcal{GH}_\Lambda(\frak{S}\times \mathbb{R}))$, sending each $f$ to $f$. In fact, the identification $\mathcal{A}_\Delta \to \mathcal{A}_{g.\Delta}$ that suits better for the formulation of the invariance under the action of ${\rm MCG}(\frak{S})$ is the map sending $f$ to $g.f$, using the action of ${\rm MCG}(\frak{S})$ on $C^\infty(\mathcal{GH}_\Lambda(\frak{S}\times \mathbb{R}))$. All the equalities appearing in the formulation of the ${\rm MCG}(\frak{S})$-invariance should be of the latter sort, taking into account the non-trivial natural action of ${\rm MCG}(\frak{S})$; so, in a sense, maybe it is a better idea to denote them by the isomorphism sign $\cong$ rather than the equality sign $=$. 

\vs

We claim that by a further investigation of the ${\rm MCG}(\frak{S})$-invariance along this line, one can obtain a version of ${\rm MCG}(\frak{S})$-equivariance of the quantization of $\mathcal{GH}_\Lambda(\frak{S}\times \mathbb{R})$ as a consequence. In particular, this will yield unitary representations of ${\rm MCG}(\frak{S})$ on Hilbert spaces, which are the principal objects of study of the present article. As mentioned in \S\ref{sec:introduction}, the main goal of the present article is to present an explicit construction of these representations of ${\rm MCG}(\frak{S})$, for it is only hinted but not fully displayed in \cite{KS}.

\section{Quantum dilogarithm functions}
\label{sec:quantum_dilogarithm}

Here we explain some special functions called the {\it quantum dilogarithm} functions, being crucial ingredients of the constructions of \cite{KS}. Here we focus mostly just on the definitions and brief partial history, and refrain from giving much further survey; for more information, see \cite{FK94, Zagier, KS, K21a}. Recall that the logarithm function with complex variable admits a power series expansion $-\log(1-z) = \sum_{m=1}^\infty \frac{z^m}{m}$ in the region $|z|<1$. Euler's dilogarithm is the function ${\rm Li}_2(z) = \sum_{m=1}^\infty \frac{z^m}{m^2}$, related to the usual logarithm by integration. Replacing one $m$ in the denominator of the $m$-th summand by the `quantum $m$', one obtains a crude version of what is called the quantum dilogarithm. It is more convenient to consider a variant ${\rm Li}_2(z;q) = \sum_{m=1}^\infty \frac{z^m}{m(1-q^m)}$, or better, its negative exponential $e^{-{\rm Li}_2(z;q)}$, which by a slight change of variables for $q$ and $z$ (namely, put $q^2$ and $-qz$ in place of $q$ and $z$) yields the function
$$
\psi^q(z) = \prod_{m=0}^\infty (1+q^{2m+1}z)^{-1}
$$
which is sometimes referred to as the {\it compact quantum dilogarithm}. It enjoys several remarkable properties, the most fundamental one being the difference equation
$$
\psi^q(q^2 z) = (1+qz) \psi^q(z),
$$
and the most interesting one being the {\it pentagon identity}: for variables $u,v$ such that $uv = q^2 vu$, 
$$
\psi^q(u) \psi^q(v) = \psi^q(v) \psi^q(q vu) \psi^q(v)
$$
holds in an appropriate sense (\cite{FK94}). 

\vs

This function $\psi^q$ has appeared in the literature for a long time in many different guises, and relatively recently in 1990's Faddeev and Kashaev started to establish some basics of it for use in mathematical physics \cite{Faddeev, FK94}. Then it was used crucially in the quantization of Teichm\"uller spaces by Kashaev \cite{Kash98} and by Chekhov and Fock \cite{CF, F}, and in the quantization of cluster varieties by Fock and Goncharov \cite{FG09}.

\vs

One way of understanding $\psi^q$ is to let $q$ be a formal variable. However, for a `physical' quantization, $q$ should be an actual number. When $q$ is a complex number satisfying $|q|<1$, the infinite product expression yields a well-defined meromorphic function on the complex plane. However, in quantization of Teichm\"uller spaces and cluster varieties \cite{Kash98, CF, F, FG09}, what was needed is the case of $|q|=1$, for which the infinite product does not make sense. Hence another version of a quantum dilogarithm function had to be used, which would share some similar properties with $\psi^q$. It is the {\em non-compact quantum dilogarithm} function $\Phi^\hbar(z)$ \cite{Faddeev, FK94}, which is a function in the complex variable $z$, associated to a nonzero real parameter $\hbar$. Assuming $\hbar>0$ at the moment for convenience, it is given by the following contour integral formula which goes back to Barnes \cite{B01} 100 years ago in the study of `double Gamma function':
\begin{align}
\label{eq:Phi_hbar}
\Phi^\hbar(z) = \exp\left( - \frac{1}{4} \int_\Omega \frac{e^{-{\rm i} pz}}{\sinh(\pi p) \sinh (\pi \hbar p)} \frac{dp}{p} \right)
\end{align}
for $z$ living in the strip $|{\rm Im}(z)| < \pi(1+\hbar)$, where $\Omega$ is the contour along the real line that avoids the origin by a small half-circle above the origin. One can show that this formula defines a non-vanishing complex analytic function on this strip, and that this function has an analytic continuation to a meromorphic function on the whole complex plane with an explicitly described sets of zeroes and poles. A fundamental property is the pair of difference equations
$$
\left\{
\begin{array}{lcr}
\Phi^\hbar(z+2\pi {\rm i} \hbar) &= & (1+e^{\pi {\rm i} \hbar} e^z) \Phi^\hbar(z), \\
\Phi^\hbar(z+2\pi {\rm i}) & = & (1+e^{\pi {\rm i} / \hbar} e^{z/\hbar}) \Phi^\hbar(z),
\end{array}
\right.
$$
which play crucial roles in solving the intertwining equations.  It also satisfies $|\Phi^\hbar(x)|=1$ for $x\in \mathbb{R}$; hence, for any self-adjoint operator ${\bf A}$ on a Hilbert space, the operator $\Phi^\hbar({\bf A})$, obtained by applying the functional calculus of ${\bf A}$ to the function $\Phi^\hbar$, is unitary. Another important property of $\Phi^\hbar$ is the pentagon identity: for self-adjoint operators ${\bf P}$ and ${\bf Q}$ on a separable Hilbert space satisfying the Heisenberg relation $[{\bf P}, {\bf Q}] = 2\pi {\rm i} \hbar \cdot {\rm id}$ (more precisely, its `Weyl-relations version' as in eq.\eqref{eq:Weyl-relations_version}), one has the following operator identity among unitary maps (\cite{FK94, FKV, G08, W00}):
$$
\Phi^\hbar({\bf P}) \Phi^\hbar({\bf Q}) = \Phi^\hbar({\bf Q}) \Phi^\hbar({\bf P}+{\bf Q}) \Phi^\hbar({\bf P}).
$$
Meanwhile, we will also need the case $\hbar<0$. The formula in eq.\eqref{eq:Phi_hbar} actually yields
$$
\Phi^{-\hbar}(z) = (\Phi^\hbar(z))^{-1},
$$
which is a way of understanding the case of a negative real parameter $\hbar$.

\vs

Comparing the difference equation of $\psi^q$ and the first of the two difference equations of $\Phi^\hbar$, keeping in mind that $q=e^{\pi {\rm i} \hbar}$, one can immediately expect that $\psi^q$ and $\Phi^\hbar$ should be closely related. To see this more clearly, let $h$ be any nonzero complex number with ${\rm Re}(h)\ge 0$. Notice that we are not using the symbol $\hbar$ here. Consider the integral expression
\begin{align}
\label{eq:Phi_h}
\Phi^h(z) = \exp\left( - \frac{1}{4} \int_{\Omega'} \frac{e^{-{\rm i} pz}}{\sinh(\pi p) \sinh (\pi h p)} \frac{dp}{p} \right),
\end{align}
which is referred to as the Barnes integral in \cite{KS}, for some contour $\Omega'$ to be explained. Let us first focus on the case when ${\rm Re}(h)>0$, for which we can take the contour $\Omega'$ to be same as before, i.e. as $\Omega$. This expression defines a non-vanishing complex analytic function on the strip $|{\rm Im}(z)|<\pi(1+{\rm Re}(h))$, and is analytically continued to a meromorphic function. When ${\rm Im}(h)>0$ is also satisfied, then by the residue computation, one can show that
\begin{align}
\label{eq:Phi_h_as_ratio}
\Phi^h(z) = \frac{\psi^{\exp(\pi {\rm i} \hbar)}(e^z)}{\psi^{\exp(-\pi {\rm i}/h)}(e^{z/h})}
\end{align}
holds, which is a straightforward exercise in complex analysis. That is, the function $\Phi^h$ is a certain ratio of two compact quantum dilogarithm functions. Notice that both quantum parameters $e^{\pi {\rm i} h}$ and $e^{-\pi {\rm i} /h}$ for the two compact quantum dilogarithms are complex numbers of modulus less than $1$, due to the assumption ${\rm Im}(h)>0$; one also notes that $h$ and $-1/h$ are related to each other by the `modular transformation'. However, notice that what we referred to as the non-compact quantum dilogarithm falls into the case when $h = \hbar$ belongs to $\mathbb{R}_{>0}$, i.e. ${\rm Im}(h)=0$. In this case, the parameters $e^{\pi {\rm i} h}$ and $e^{-\pi {\rm i} /h}$ are of modulus $1$, hence the associated compact quantum dilogarithm functions do not make sense, therefore nor does the equality in eq.\eqref{eq:Phi_h_as_ratio}. In a sense, the non-compact quantum dilogarithm can be thought of as being a limit of the ratio of compact quantum dilogarithms as in the right hand side of eq.\eqref{eq:Phi_h_as_ratio} as ${\rm Im}(h)\to 0$. Anyways, the function $\Phi^h(z)$ defined in eq.\eqref{eq:Phi_h} can also be referred to as a quantum dilogarithm function as well.

\vs

The case ${\rm Re}(h)>0$ has been studied in the literature, but for the case ${\rm Re}(h)=0$, i.e. when $h = {\rm i}\hbar \in {\rm i} \mathbb{R}$ is a pure imaginary number, the integral expression in eq.\eqref{eq:Phi_h} seems to have been somewhat elusive. In \cite{KS} we undertook the task of studying this case. The contour integral in eq.\eqref{eq:Phi_h} does not work if we use the previous $\Omega$ as the contour $\Omega'$, because in case ${\rm Re}(h)=0$ the integrand has poles along real and imaginary axes, hence also on $\Omega$. Thus we used a `slanted' contour $\Omega' = e^{{\rm i} \theta} \Omega$ for some angle $\theta \in \mathbb{R}$ of small absolute value, which is taken to be a negative angle when ${\rm Im}(h) \ge 0$ and a positive angle when ${\rm Im}(h) \le 0$. We then still get a well-defined meromorphic function $\Phi^h = \Phi^{{\rm i} \hbar}$, and when ${\rm Im}(h)>0$ the ratio expression in eq.\eqref{eq:Phi_h_as_ratio} still holds; in particular, when ${\rm Im}(h)>0$ and ${\rm Re}(h)=0$, i.e. $h = {\rm i} \hbar \in {\rm i}\mathbb{R}_{>0}$, we have:
$$
\Phi^{{\rm i}\hbar}(z) = \frac{\psi^{\exp(-\pi \hbar)}(e^z)}{\psi^{\exp(-\pi/\hbar)}(e^{z/({\rm i} \hbar)})}.
$$
Notice that the quantum parameters $e^{-\pi\hbar}$ and $e^{-\pi/\hbar}$ for the two compact quantum dilogarithms in this case are positive real numbers less than $1$. Similar arguments show that for ${\rm Im}(h)<0$, i.e. when $h = - {\rm i}\hbar$ with $\hbar>0$, the function $\Phi^{-{\rm i}\hbar}(z)$ can be computed to be 
$$
\Phi^{-{\rm i}\hbar}(z) = \frac{\psi^{\exp(-\pi/\hbar)}(e^{-z/({\rm i}\hbar)})}{\psi^{\exp(-\pi\hbar)}(e^z)} = \ol{\Phi^{{\rm i}\hbar}(\ol{z})}^{-1}.
$$
One can view each of $\Phi^{{\rm i}\hbar}$ and $\Phi^{-{\rm i} \hbar}$ as a `modular double' version of a compact quantum dilogarithm function $\psi^{\bf q}$, with a real quantum parameter ${\bf q} = e^{-\pi \hbar}$. As mentioned, although these functions $\Phi^{{\rm i}\hbar}$ and $\Phi^{-{\rm i} \hbar}$ are genuine ratios of certain compact quantum dilogarithm functions, they do not seem to have been considered before in the literature.

\vs

What is actually introduced and used in \cite{KS} for the quantization of $\mathcal{GH}_\Lambda(\frak{S}\times\mathbb{R})$ is a two-real-variable function $F^\hbar_\Lambda : \mathbb{R}^2 \to \mathbb{C}$, depending on the cosmological constant $\Lambda \in \{-1,0,1\}$, defined by
\begin{align}
\label{eq:F_hbar_Lambda}
F^\hbar_\Lambda(x,y) = \left\{
{\renewcommand{\arraystretch}{1.4} \begin{array}{ll}
\Phi^\hbar(x+ \hbar y) \Phi^{-\hbar}(x- \hbar y) & \mbox{if $\Lambda=-1$}, \\
\Phi^{{\rm i} \hbar}(x+{\rm i} \hbar y) \Phi^{-{\rm i} \hbar}(x - {\rm i} \hbar y) & \mbox{if $\Lambda=1$}, \\
(1+e^x)^{y/(\pi {\rm i})} = \exp\left( \frac{y}{\pi {\rm i}}\log(1+e^x) \right) & \mbox{if $\Lambda=0$}.
\end{array}}
\right.
\end{align}
The three formulas may not look so uniform at first glance. We note that in fact they all come from the heuristic expression
$$
\Phi^{\ell \hbar}(x+\ell \hbar y) \Phi^{-\ell \hbar}(x-\ell \hbar y)
$$
in terms of the ring $\mathbb{R}_\Lambda = \mathbb{R}[\ell]/(\ell^2 = - \Lambda)$ of generalized complex numbers, or in terms of a suitable complexification of $\mathbb{R}_\Lambda$. Note also that the second factor $\Phi^{-\ell\hbar}(x-\ell\hbar y)$ is obtained by applying the $\mathbb{R}_\Lambda$-conjugation on the first factor $\Phi^{\ell\hbar}(x+\ell\hbar y)$, so the first and the second factors can be thought of as $\mathbb{R}_\Lambda$-holomorphic and $\mathbb{R}_\Lambda$-anti-holomorphic parts. But this is indeed just at a heuristic level, as the {\it $\mathbb{R}_\Lambda$-version of the quantum dilogarithm} $\Phi^{\ell\hbar}$ cannot be defined using the contour integral on the nose. Anyways, the function $F^\hbar_0$ for $\Lambda=0$ `comes from' an $\mathbb{R}_0$-version of the quantum dilogarithm, hence can be viewed as a degenerate, or flat, version of the quantum dilogarithm.

\vs

Unlike $\Phi^{\pm \ell\hbar}$, the defining expressions in the right hand side of eq.\eqref{eq:F_hbar_Lambda} for the function $F^\hbar_\Lambda$ do make good senses. One useful property of the function $F^\hbar_\Lambda$ we shall use is its `unitarity':
\begin{align}
\label{eq:F_hbar_Lambda_unitarity}
|F^\hbar_\Lambda(x,y)| =1, \qquad \forall x,y\in \mathbb{R}.
\end{align}
If ${\bf A}$ and ${\bf B}$ are self-adjoint operators on a Hilbert space that strongly commute with each other, i.e. if they satisfy the Weyl-relations version (eq.\eqref{eq:Weyl-relations_version}) of the Heisenberg relation $[{\bf A}, {\bf B}]=0$, then one can apply the two-variable functional calculus for ${\bf A}$ and ${\bf B}$ to this function $F^\hbar_\Lambda$ to obtain an operator $F^\hbar_\Lambda({\bf A}, {\bf B})$. For the multi-variable functional calculus for a strongly commuting collection of self-adjoint operators, see \cite{Sch}. The unitarity as in eq.\eqref{eq:F_hbar_Lambda_unitarity} implies that $F^\hbar_\Lambda({\bf A}, {\bf B})$ is unitary. Another important property of $F^\hbar_\Lambda$ is the pentagon operator identity, as shown in \cite[Prop.6.16]{KS}. It says that, if ${\bf x}_1$, ${\bf y}_1$, ${\bf x}_2$ and ${\bf y}_2$ are self-adjoint operators satisfying the (Weyl-relations versions of the) Heisenberg relations $[{\bf x}_1, {\bf x}_2]=[{\bf y}_1,{\bf y}_2]=[{\bf x}_1,{\bf y}_1]=[{\bf x}_2,{\bf y}_2]=0$, $[{\bf x}_1, {\bf y}_2] = [{\bf y}_1, {\bf x}_2] = \pi {\rm i} \cdot {\rm id}$, then 
\begin{align}
\label{eq:F_hbar_Lambda_pentagon}
F^\hbar_\Lambda({\bf x}_1,{\bf y}_1) F^\hbar_\Lambda({\bf x}_2,{\bf y}_2)
= F^\hbar_\Lambda({\bf x}_2,{\bf y}_2)
F^\hbar_\Lambda({\bf x}_1+{\bf x}_2,{\bf y}_1+{\bf y}_2) 
F^\hbar_\Lambda({\bf x}_1,{\bf y}_1) 
\end{align}
holds as an equality of unitary operators. For other properties of $F^\hbar_\Lambda$, see \cite{KS}.

\section{The unitary interwiners}
\label{sec:the_unitary_intertwiners}

Among the aspects of the quantization program of \cite{KS} as spelled out in \S\ref{sec:overview}, in the present article we focus only on the system of Hilbert spaces $\mathscr{H}_\Delta$ of quantum states associated to ideal triangulations $\Delta$, together with the unitary maps ${\bf K}^\hbar_{\Delta,\Delta'}$ which identify these Hilbert spaces with each other. For the remaining aspects of the quantization program, see \cite{KS}.

\subsection{The Hilbert space and basic operators}

Fix any ideal triangulation $\Delta$. Let $\mathscr{H}_\Delta$ be a separable complex Hilbert space, and let ${\bf x}_i$ and ${\bf y}_i$ ($i\in \Delta$) be self-adjoint operators on $\mathscr{H}_\Delta$ satisfying the Heisenberg relations
\begin{align}
\label{eq:x_i_y_i_Heisenberg}
[{\bf x}_i, {\bf x}_j]=0, \qquad
[{\bf y}_i, {\bf y}_j]=0, \qquad
[{\bf x}_i, {\bf y}_j]=\pi {\rm i} \varepsilon_{ij} \cdot {\rm id}, \qquad \forall i,j\in \Delta,
\end{align}
One might wonder why there is no quantum parameter $\hbar$ in these Heisenberg relations. Indeed, what should be quantum counterparts of the coordinate functions $x_i$ and $y_i$ of $\mathcal{GH}_\Lambda(\frak{S}\times \mathbb{R})$ are the operators ${\bf x}_i$ and $\hbar {\bf y}_i$, as we have $[{\bf x}_i, \hbar {\bf y}_j] = \pi {\rm i} \hbar \varepsilon_{ij} \cdot {\rm id}$. In \cite{KS} we decided to take out $\hbar$, which made it easier to `$\mathbb{R}_\Lambda$-complexify' $\hbar\in \mathbb{R}$ to $\ell \hbar \in \mathbb{R}_\Lambda$.

\vs

When we say that two self-adjoint operators ${\bf A}$ and ${\bf B}$ satisfy a Heisenberg relation
$$
[{\bf A}, {\bf B}] = {\rm i} c \cdot {\rm id}
$$
for some $c\in \mathbb{R}$, unless otherwise stated we mean that they satisfy the {\it Weyl-relations} version of it, i.e. that the family of equations
\begin{align}
\label{eq:Weyl-relations_version}
e^{{\rm i}\alpha {\bf A}} e^{{\rm i}\beta {\bf B}}
= e^{-{\rm i}c\alpha\beta} e^{{\rm i}\beta {\bf B}} e^{{\rm i}\alpha {\bf A}}
\end{align}
of unitary operators are satisfied for all $\alpha,\beta \in \mathbb{R}$. This is what makes the solution of the Heisenberg relations more rigid, and is required as a hypothesis for a (generalized) Stone-von Neumann theorem. Let us also consider the quantum constraint equations
\begin{align}
\label{eq:quantum_constraint_equations}
\sum_{i \in \Delta} v_{i,p} {\bf x}_i = 0 \qquad\mbox{and}\qquad
\sum_{i \in \Delta} v_{i,p} {\bf y}_i = 0, \qquad \forall\mbox{puncture $p$}.
\end{align}
One can observe that the Stone-von Neumann theorem \cite{vN} \cite[\S14]{Hall} implies that the collection of self-adjoint operators $\{ {\bf x}_i, {\bf y}_i\}_{i\in \Delta}$ on a separable complex Hilbert space $\mathscr{H}$ satisfying eq.\eqref{eq:x_i_y_i_Heisenberg} and eq.\eqref{eq:quantum_constraint_equations} is uniquely determined up to unitary equivalence, if we also stipulate that these operators should act in an {\it irreducible} manner, meaning that the only closed subspaces of $\mathscr{H}$ invariant under $e^{{\rm i} \alpha {\bf x}_i}$ and $e^{{\rm i} \alpha {\bf y}_i}$ for all $\alpha \in \mathbb{R}$ and $i$ are $0$ and $\mathscr{H}$ itself.

\vs

As just mentioned, by using the existence and uniqueness result of the Stone-von Neumann theorem, it is enough to just write the equations eq.\eqref{eq:x_i_y_i_Heisenberg} and eq.\eqref{eq:quantum_constraint_equations}, without writing down the actual Hilbert space and the operators explicitly; see \cite{K21b} \cite[\S14]{Hall}. However, instead of doing so, we shall explicitly present a particular collection of operators as a solution. For convenience of our presentation, we will at the moment forget the constraint equations in eq.\eqref{eq:quantum_constraint_equations} and work with a certain solution to just the Heisenberg relations in eq.\eqref{eq:x_i_y_i_Heisenberg}, based on which we will build our unitary intertwiners. Then afterwards in \S\ref{subsec:irreducible_representations}, we will give a recipe how to actually obtain a solution satisfying both eq.\eqref{eq:x_i_y_i_Heisenberg} and eq.\eqref{eq:quantum_constraint_equations}, and the corresponding intertwiners.

\vs

The Hilbert space is set to be
$$
\mathscr{H}_\Delta := L^2(\mathbb{R}^\Delta, \underset{i \in \Delta}{\textstyle \bigwedge} dt_i),
$$
with respect to the standard Lebesgue measure on $\mathbb{R}^\Delta = \mathbb{R}^{|\Delta|} = \mathbb{R}^{6g-6+3n}$. Consider the position and the momentum operators
$$
t_i \qquad\mbox{and}\qquad {\rm i} \frac{\partial}{\partial t_i}, \qquad \mbox{for } ~ i \in \Delta,
$$
which are self-adjoint operators on $\mathscr{H}_\Delta$ satisfying the Heisenberg relations
$$
[t_i, t_j] = 0, \qquad
[{\rm i} \frac{\partial}{\partial t_i}, {\rm i} \frac{\partial}{\partial t_j}]=0,\qquad
[t_i, {\rm i} \frac{\partial}{\partial t_j}] = - {\rm i} \delta_{i,j} \cdot {\rm id},
$$
where $\delta_{i,j}$ is the Kronecker delta. It is well known that any $\mathbb{R}$-linear combinations of the position and the momentum operators yield self-adjoint operators, and any two of them satisfy the Heisenberg relations that can be deduced from the bilinearity of the commutator bracket; see \cite{Hall}. Using this fact, we consider the operators
\begin{align}
\label{eq:x_i_y_i_operators}
{\bf x}_i = - \pi {\rm i} \frac{\partial}{\partial t_i}, \qquad
{\bf y}_i = \underset{j \in \Delta}{\textstyle \sum} \varepsilon_{ij} t_j, \qquad \forall i \in \Delta,
\end{align}
which are self-adjoint operators on $\mathscr{H}_\Delta$ satisfying the Heisenberg relations in eq.\eqref{eq:x_i_y_i_Heisenberg}.

\subsection{The intertwiner for a flip at an arc}
\label{subsec:flip_intertwiner}

Suppose that $\Delta$ and $\Delta'$ are ideal triangulations of $\frak{S}$ that are related to each other by {\it flip} at an arc $k$, that is, they differ exactly at one arc denoted by $k$. In this case, there is a natural set bijection between $\Delta$ and $\Delta'$, allowing us to enumerate the members of these two triangulations by the same symbols like $i$, $j$ and $k$. We aim to describe the solution to the unitary intertwiner
$$
{\bf K}^\hbar_{\Delta,\Delta'} : \mathscr{H}_{\Delta'} \to \mathscr{H}_\Delta
$$
constructed in \cite{KS}. Since $\Delta'$ is obtained from $\Delta$ by flip at $k$, we may write $\Delta' = \mu_k(\Delta)$ using the symbol $\mu_k$ which hints at the connection to the `cluster mutations' \cite{FZ, FG06, FG09}, and denote this intertwiner as
$$
{\bf K}^\hbar_k = {\bf K}^\hbar_{\Delta,\Delta'}= {\bf K}^\hbar_{\Delta,\mu_k(\Delta)} : \mathscr{H}_{\mu_k(\Delta)} = \mathscr{H}_{\Delta'} \to\mathscr{H}_\Delta.
$$
One has to be careful when dealing with the symbol ${\bf K}^\hbar_k$ because it is ambiguous, as it does not indicate the underlying ideal triangulations. 

\vs

Following \cite{FG09}, the intertwiner ${\bf K}^\hbar_k = {\bf K}^\hbar_{\Delta,\Delta'}$ for a flip $\mu_k$ at $k$ is obtained as the composition
$$
\xymatrix{
{\bf K}^\hbar_k = {\bf K}^{\hbar\sharp}_k \circ {\bf K}'_k ~:~ \mathscr{H}_{\Delta'} \ar[r]^-{{\bf K}'_k} & \mathscr{H}_\Delta \ar[r]^-{{\bf K}^{\hbar\sharp}_k} & \mathscr{H}_\Delta
}
$$
of two unitary operators ${\bf K}'_k$ and ${\bf K}^{\hbar\sharp}_k$, where
$$
{\bf K}'_k = {\bf K}'_{\Delta,\Delta'} : \mathscr{H}_{\Delta'} = L^2(\mathbb{R}^{\Delta'}, \underset{i \in \Delta'}{\textstyle \bigwedge} dt_i') \to L^2(\mathbb{R}^\Delta, \underset{i \in \Delta}{\textstyle \bigwedge} dt_i) = \mathscr{H}_\Delta
$$
is called the `monomial transformation part' and
$$
{\bf K}^{\hbar\sharp}_k = {\bf K}^{\hbar\sharp}_{\Delta,\Delta'} : \mathscr{H}_\Delta = L^2(\mathbb{R}^\Delta, \underset{i \in \Delta}{\textstyle \bigwedge} dt_i) \to L^2(\mathbb{R}^\Delta, \underset{i \in \Delta}{\textstyle \bigwedge} dt_i) = \mathscr{H}_\Delta
$$
is called the `automorphism part'. 

\vs

The monomial transformation part ${\bf K}'_k : L^2(\mathbb{R}^{\Delta'}, \underset{i \in \Delta'}{\textstyle \bigwedge} dt_i') \to L^2(\mathbb{R}^\Delta, \underset{i \in \Delta}{\textstyle \bigwedge} dt_i)$ is the operator induced by the map
$$
\chi_k : \mathbb{R}^{\Delta} \to \mathbb{R}^{\Delta'}
$$
between the Euclidean spaces $\mathbb{R}^{\Delta}$ and $\mathbb{R}^{\Delta'}$ defined as follows: the pullback of each coordinate function of the codomain is given by
$$
\chi_k^* t_i' = \left\{
\begin{array}{ll}
t_i & \mbox{if $i\neq k$}, \\
-t_k + \sum_{j\in \Delta} [-\varepsilon_{kj}]_+ \, t_j & \mbox{if $i=k$},
\end{array}
\right.
$$
where $\varepsilon=(\varepsilon_{ij})_{i,j\in \Delta}$ is the {\it exchange matrix} for the triangulation $\Delta$ defined as
\begin{align*}
\varepsilon_{ij} & = a_{ij} - a_{ji}, \\
a_{ij} & = \mbox{the number of corners of ideal triangles of $\Delta$} \\
& \hspace{4,8mm} \mbox{delimited by $i$ at the left and $j$ at the right},
\end{align*}
and $[a]_+$ is the positive part of a real number $a$, i.e.
$$
[a]_+ = \frac{a+|a|}{2}.
$$
In \cite{K21a}, the author studied the operators on $L^2(\mathbb{R}^N)$ induced by invertible linear maps on $\mathbb{R}^N$, precisely to the extent we need in this setting; in particular, as observed in \cite{K21a}, the operator ${\bf K}'_k$ is unitary.

\vs

The automorphism part ${\bf K}^{\hbar\sharp}_k : L^2(\mathbb{R}^{\Delta}, \underset{i \in \Delta}{\textstyle \bigwedge} dt_i) \to L^2(\mathbb{R}^\Delta, \underset{i \in \Delta}{\textstyle \bigwedge} dt_i)$ is the operator given by the following formula, depending on the cosmological constant $\Lambda$:
\begin{align}
\label{eq:K_sharp_formula}
{\bf K}^{\hbar \sharp}_k = F^\hbar_\Lambda({\bf x}_k, {\bf y}_k) = \left\{
{\renewcommand{\arraystretch}{1.4} \begin{array}{ll}
\Phi^\hbar({\bf x}_k + \hbar {\bf y}_k) \, \Phi^{-\hbar}({\bf x}_k - \hbar {\bf y}_k) & \mbox{if $\Lambda=-1$}, \\
\Phi^{{\rm i}\hbar}({\bf x}_k + {\rm i}\hbar {\bf y}_k) \Phi^{-{\rm i}\hbar}({\bf x}_k - {\rm i} \hbar {\bf y}_k) & \mbox{if $\Lambda=1$}, \\
(1+e^{ {\bf x}_k})^{{\bf y}_k/(\pi {\rm i})} & \mbox{if $\Lambda=0$},
\end{array}}
\right.
\end{align}
where the function $F^\hbar_\Lambda : \mathbb{R}^2 \to \mathbb{C}$ is as defined in eq.\eqref{eq:F_hbar_Lambda} in \S\ref{sec:quantum_dilogarithm}. Here $\Phi^{\pm \hbar}$ and $\Phi^{\pm {\rm i}\hbar}$ are some versions of the quantum dilogarithm functions described in \S\ref{sec:quantum_dilogarithm}. The formula $F^\hbar_\Lambda({\bf x}_k, {\bf y}_k)$, i.e. each of the three lines in the right of eq.\eqref{eq:K_sharp_formula}, is defined by the two-variable functional calculus for the self-adjoint operators ${\bf x}_k$ and ${\bf y}_k$ on $\mathscr{H}_\Delta = L^2(\mathbb{R}^{\Delta}, \underset{i \in \Delta}{\textstyle \bigwedge} dt_i)$ in eq.\eqref{eq:x_i_y_i_operators}, which makes sense because ${\bf x}_k$ and ${\bf y}_k$ strongly commute with each other; indeed, recall that they satisfy the Weyl-relations version of the Heisenberg relation $[{\bf x}_k, {\bf y}_k]=0$, as in eq.\eqref{eq:x_i_y_i_Heisenberg}.  From eq.\eqref{eq:F_hbar_Lambda_unitarity}, one observes that ${\bf K}^{\hbar \sharp}_k$ is unitary.

\vs

As said in \S\ref{sec:quantum_dilogarithm}, we have $\Phi^{-\hbar}(z) = (\Phi^\hbar(z))^{-1}$, which entails that the above described intertwiner ${\bf K}^\hbar_k = {\bf K}^{\hbar \sharp}_k \circ {\bf K}'_k$ for the case $\Lambda=-1$ coincides with Fock and Goncharov's intertwiner associated to a cluster mutation $\mu_k$ \cite{FG09}, for the quantization of their `cluster symplectic double variety', or `cluster $\mathscr{D}$-variety' \cite{FG09}. The construction of \cite{FG09} applies to {\em any} skew-symmetric integer matrix $\varepsilon$. If we restrict ourselves to the case when the exchange matrix $\varepsilon$ comes from an ideal triangulation of a punctured surface, what Fock and Goncharov quantized is a certain double of the `enhanced Teichm\"uller space' $\mathcal{T}^+(\frak{S})$ \cite{BL, F}, which in a sense can be viewed as the cotangent bundle of $\mathcal{T}^+(\frak{S})$. %Here, the enhanced Teichm\"uller space $\mathcal{T}^+(\frak{S})$ is what parametrizes a complete hyperbolic metric on $\frak{S}$ up to isotopy, together with the choice of an orientation (a choice in the set $\{+,-\}$) per each funnel end \cite{BL}. 
Anyhow, we note that the functions used in this case of $\Lambda=-1$ are the non-compact quantum dilogarithm functions $\Phi^{\pm \hbar}$ of Faddeev and Kashaev \cite{Faddeev, FK94, B01}, which are also what are used in the original works of quantum Teichm\"uller theory by Kashaev \cite{Kash98} and by Chekhov and Fock \cite{CF, F} for the construction of unitary intertwiners; in these works, the flip intertwiner involves just one copy of $\Phi^\hbar$, instead of both $\Phi^\hbar$ and $\Phi^{-\hbar} = (\Phi^\hbar)^{-1}$.

\subsection{The intertwiner for an arbitrary change of ideal triangulations}

Now let $\Delta$ and $\Delta'$ be any two ideal triangulations of $\frak{S}$. It is well known that they can be connected to each other by a sequence of flips; see Labardini-Fragoso \cite{LF} for such a statement when we restrict only to the ideal triangulations without self-folded triangles. So one can find a sequence of ideal triangulations $\Delta_0,\Delta_1,\ldots,\Delta_m$, such that $\Delta_0=\Delta$, $\Delta_m = \Delta'$, and $\forall i=1,\ldots,m$, $\Delta_{i-1}$ and $\Delta_i$ are related by a flip, say $\Delta_i = \mu_{k_i}(\Delta_{i-1})$ for some arc $k_i$. We then define the intertwiner ${\bf K}^\hbar_{\Delta,\Delta'}$ associated to the change of ideal triangulations $\Delta \leadsto \Delta'$ as the composition of the intertwiners for flips in the above sequence. That is, we let
\begin{align}
\label{eq:intertwiner_as_composition_of_sequence}
{\bf K}^\hbar_{\Delta,\Delta'} := {\bf K}^\hbar_{\Delta_0,\Delta_1} \circ {\bf K}^\hbar_{\Delta_1,\Delta_2} \circ \cdots \circ {\bf K}^\hbar_{\Delta_{m-1},\Delta_m}  = {\bf K}^\hbar_{k_1} \circ {\bf K}^\hbar_{k_2} \circ \cdots \circ {\bf K}^\hbar_{k_m} : \mathscr{H}_{\Delta'} \to \mathscr{H}_\Delta
\end{align}

\vs

As the sequence of flips connecting $\Delta$ and $\Delta'$ is not unique, one has to make sure that the resulting operator ${\bf K}^\hbar_{\Delta,\Delta'}$ does not depend on the choice of such a sequence. This well-definedness is guaranteed if we show that the flip operators ${\bf K}^\hbar_k$ satisfy all possible relations satisfied by the flips of ideal triangulations. The first step is to classify all relations of flips of triangulations. The following result is well known; see \cite[Thm.3.10]{FST} and references therein.
\begin{proposition}
\label{prop:relations_of_flips}
Any relation satisfied by flips of ideal triangulations without self-folded triangles is generated by the following ones:
\begin{enumerate}[label=\rm (\arabic*)]
\item (twice-flip identity) $\mu_k \circ \mu_k = {\rm id}$ when applied to any $\Delta$;

\item (quadrilateral identity) $\mu_i \circ \mu_j \circ  \mu_i \circ \mu_j = {\rm id}$ when applied to any $\Delta$ with $\varepsilon_{ij}=0$;

\item (pentagon identity) $P_{(ij)} \circ \mu_i \circ \mu_j \circ \mu_i \circ \mu_j \circ \mu_i = {\rm id}$ when applied to any $\Delta$ with $\varepsilon_{ij}=\pm 1$, where $P_{(ij)}$ means re-labeling $i\leftrightarrow j$.
\end{enumerate}
\end{proposition}
The last relation actually holds without the re-labeling process. However, we find it useful to keep track of the labels of ideal arcs, in which case we are dealing in fact with labeled ideal triangulations in effect. We fix an index set $I$, say $\{1,2,\ldots,6g-6+3n\}$, and by a labeled ideal triangulation we mean an ideal triangulation $\Delta$ together with a set bijection $\Delta \to I$, which is the labeling map. By a slight abuse of notation, we still denote just by the same symbol $\Delta$ a labeled ideal triangulation. Each element of $I$ corresponds to an ideal arc of a labeled ideal triangulation; so the symbols $i$, $j$ and $k$ which appeared so far can be understood as elements of $I$. The Hilbert space $\mathscr{H}_\Delta = L^2(\mathbb{R}^\Delta, \bigwedge_{i\in \Delta} dt_i)$ can now be written as
$$
\mathscr{H}_\Delta = L^2(\mathbb{R}^I, \underset{i \in I}{\textstyle \bigwedge} dt_i).
$$
When we deal with several labeled ideal triangulations $\Delta$ at the same time, we will denote the variables $t_i$ for $\mathscr{H}_\Delta$ differently for different $\Delta$, if necessary.

\vs

Now, in addition to the flip $\mu_k$ at an arc $k\in I$, we add one more kind of elementary transformation of labeled ideal triangulations, which is the label permutation $P_\sigma$ associated to any permutation $\sigma$ of the index set $I$. We write $\Delta' = P_\sigma(\Delta)$ if $\Delta$ and $\Delta'$ are the same ideal triangulations with their labeling maps related as follows: if the constituent ideal arc is labeled by $i$ in $\Delta$ then it is labeled by $\sigma(i)$ in $\Delta'$. To this change of ideal triangulations $\Delta \leadsto \Delta' = P_\sigma(\Delta)$ we associate the index-permutation sintertwiner 
$$
{\bf P}_\sigma = {\bf K}^\hbar_{\Delta,P_\sigma(\Delta)} = {\bf K}^\hbar_{\Delta,\Delta'} ~:~ \mathscr{H}_{\Delta'} = L^2(\mathbb{R}^I, \underset{i \in I}{\textstyle \bigwedge} dt_i') \to L^2(\mathbb{R}^I, \underset{i \in I}{\textstyle \bigwedge} dt_i) = \mathscr{H}_\Delta
$$
defined as the operator induced by the re-labeling map $L_\sigma : \mathbb{R}^I = \mathbb{R}^\Delta \to \mathbb{R}^{\Delta'} = \mathbb{R}^I$, whose pullback action on the coordinate functions is given by
$$
L_\sigma^* t'_{\sigma(i)} = t_i, \qquad \forall i \in I.
$$
It's easy to see that ${\bf P}_\sigma$ is unitary.

\vs

Now that we decided to work with labeled ideal triangulations, eq.\eqref{eq:intertwiner_as_composition_of_sequence} has to be modified because now a change of ideal triangulation $\Delta \leadsto \Delta'$ is expressed as a sequence of mutations and label permutations; hence ${\bf K}^\hbar_{\Delta,\Delta'}$ should be defined as composition of corresponding flip intertwiners ${\bf K}^\hbar_k$ and index permutation operators ${\bf P}_\sigma$. To ensure the well-definedness of ${\bf K}^\hbar_{\Delta,\Delta'}$, we now need to classify all relations satisfied by flips and label permutations.
\begin{proposition}
\label{prop:relations_permutations}
Any relation satisfied by flips and label permutations of labeled ideal triangulations without self-folded triangles is generated by the relations in Prop.\ref{prop:relations_of_flips}, together with the following:
$$
P_{\rm id} = {\rm id}, \qquad
P_{(\sigma\circ \gamma)^{-1}} \circ P_\sigma \circ P_\gamma = {\rm id}, \qquad
P_\sigma \circ \mu_i \circ P_{\sigma^{-1}} \circ \mu_{\sigma(i)}={\rm id},
$$
applied to any labeled ideal triangulation, for all permutations $\sigma$ and $\gamma$ of $I$ and for all $i\in I$.
\end{proposition}
This proposition would follow from Prop.\ref{prop:relations_of_flips} by an argument analogous to the one in \cite[\S4.1]{K16b}, using a lemma by Bakalov and Kirillov \cite{BK}.

\vs

As desired, all relations satisfied by flips and index permutations are indeed shown to be satisfied by the corresponding intertwining operators. Note that due to the contravariant nature of the construction of the intertwining operators, we need to reverse the order of composition in the relations when translating to relations of operators. For example, the operator identity corresponding to the pentagon identity in Prop.\ref{prop:relations_of_flips}(3) reads
\begin{align}
\label{eq:pentagon_operator_identity}
{\bf K}^\hbar_i \circ {\bf K}^\hbar_j \circ {\bf K}^\hbar_i \circ {\bf K}^\hbar_j \circ {\bf K}^\hbar_i \circ {\bf P}_{(ij)} = {\rm id} ~:~ \mathscr{H}_\Delta \to \mathscr{H}_\Delta
\end{align}
which holds when the `initial' triangulation $\Delta$ satisfies $\varepsilon_{ij}=\pm1$; the core of a proof of this identity is eq.\eqref{eq:F_hbar_Lambda_pentagon}, the pentagon operator identity of the function $F^\hbar_\Lambda$ which is used in the definition of the flip intertwiners. 
\begin{theorem}[consistency equations of intertwiners \cite{KS}]
\label{thm:consistency_equations_of_intertwiners}
The intertwining operators ${\bf K}^\hbar_k$ and ${\bf P}_\sigma$ corresponding to the flip $\mu_k$ and label permutation $P_\sigma$ of labeled ideal triangulations satisfy all relations appearing in Prop.\ref{prop:relations_of_flips}--\ref{prop:relations_permutations}, with the order of composition reversed.
\end{theorem}
This implies the well-definedness of the intertwiner ${\bf K}^\hbar_{\Delta,\Delta'}$ associated to an arbitrary change of labeled ideal triangulations $\Delta \leadsto \Delta'$; that is, the operator ${\bf K}^\hbar_{\Delta,\Delta'}$ does not depend on the choice of a sequence of flips and label permutations connecting $\Delta$ and $\Delta'$. Moreover, this result also implies the consistency equations of the resulting intertwiners. Namely,
\begin{align}
\label{eq:K_consistency}
{\bf K}^\hbar_{\Delta,\Delta'} \circ {\bf K}^\hbar_{\Delta',\Delta''} = {\bf K}^\hbar_{\Delta,\Delta''}
\end{align}
holds for any triple of labeled ideal triangulations $\Delta$, $\Delta'$ and $\Delta''$. 

\vs

This finishes the description of the system of Hilbert spaces $\mathscr{H}_\Delta$ of quantum states associated to labeled ideal triangulations $\Delta$, together with a consistent system of unitary intertwiners ${\bf K}^\hbar_{\Delta,\Delta'}$ among these Hilbert spaces, for the quantization program for the moduli space $\mathcal{GH}_\Lambda(\frak{S}\times \mathbb{R})$ established in \cite{KS}.

\section{Mapping class group representations}
\label{sec:MCG_representations}

We explain how the unitary intertwiners can be used to construct representations of the mapping class group ${\rm MCG}(\frak{S})$ of $\frak{S}$. As a matter of fact, perhaps a better viewpoint is to interpret the intertwiners as forming a representation of a groupoid rather than that of a group, as pointed out to the author by Dennis Sullivan at Igor Frenkel's 70th birthday conference. Indeed, in \cite{FG09} Fock and Goncharov first formulate their quantization result in terms of a groupoid, and as a consequence obtains a representation of a group. Here we do likewise. We note that the constructions presented in this section are not explicitly written down in \cite{KS}.

\subsection{The Ptolemy groupoid, and the Ptolemy class groupoid}
\label{subsec:Ptolemy_groupoid}

Following \cite{P87} with some modification, one may define the {\it Ptolemy groupoid} ${\rm Pt}(\frak{S})$ of the surface $\frak{S}$ as the category whose set of objects ${\rm Ob}({\rm Pt}(\frak{S}))$ is the set of all labeled ideal triangulations $\Delta$ without self-folded triangles, such that for each pair of objects $\Delta$ and $\Delta'$ there is exactly one morphism from $\Delta$ to $\Delta'$, denoted by $(\Delta \leadsto \Delta')$. Here, an ideal triangulation $\Delta$ is considered up to simultaneous isotopy of its constituent ideal arcs.

\vs

One can formulate the results reviewed in \S\ref{sec:the_unitary_intertwiners} as a contravariant functor
$$
{\rm Pt}(\frak{S}) \to {\rm Hilb},
$$
where ${\rm Hilb}$ denotes the category of Hilbert spaces whose morphisms are unitary maps. The object $\Delta$ of ${\rm Pt}(\frak{S})$ goes to the Hilbert space $\mathscr{H}_\Delta$, and the morphism $(\Delta\leadsto \Delta')$ of ${\rm Pt}(\frak{S})$ goes to the unitary map ${\bf K}^\hbar_{\Delta,\Delta'} : \mathscr{H}_{\Delta'} \to \mathscr{H}_\Delta$.

\vs

An important property of this functor is the invariance under the action of the mapping class group ${\rm MCG}(\frak{S}) = {\rm Diff}_+(\frak{S})/{\rm Diff}(\frak{S})_0$, which we now describe. First, we note that ${\rm MCG}(\frak{S})$ naturally acts on the Ptolemy groupoid ${\rm Pt}(\frak{S})$. The element $\h \in {\rm MCG}(\frak{S})$ naturally sends the object $\Delta$ to $\h.\Delta$ because $\h$ is an isotopy class of homeomorphisms or diffeomorphisms from $\frak{S}$ to itself, and accordingly sends the morphism $(\Delta\leadsto \Delta')$ to $\h.(\Delta\leadsto \Delta') := (\h.\Delta \leadsto \h.\Delta')$. Recall that an object $\Delta$ of ${\rm Pt}(\frak{S})$ is not just an ideal triangulation, but is equipped with a labeling map $\Delta \to I$. The object $\h.\Delta$ then is given the natural labeling map given by the composition $\h.\Delta \to \Delta \to I$; that is, if $\alpha$ is an ideal arc of $\Delta$ labeled by $i \in I$, then the ideal arc $\h.\alpha$ of $\h.\Delta$ is labeled by $i$. This induces a natural identification between the Hilbert spaces
$$
\mathscr{H}_{\h.\Delta} = L^2(\mathbb{R}^I, \underset{i \in I}{\textstyle \bigwedge} dt_i) \to L^2(\mathbb{R}^I, \underset{i \in I}{\textstyle \bigwedge} dt_i) = \mathscr{H}_\Delta
$$
which is the identity map $L^2(\mathbb{R}^I, {\textstyle \bigwedge}_{i \in I} dt_i) \to L^2(\mathbb{R}^I, {\textstyle \bigwedge}_{i \in I} dt_i)$. Via this identification, one can regard
$$
\mathscr{H}_{\h.\Delta} \equiv \mathscr{H}_\Delta.
$$
Now, in case $(\Delta\leadsto \Delta')$ is a flip $\mu_k$ at $k\in I$, recall that the construction of the corresponding intertwiner ${\bf K}^\hbar_{\Delta,\Delta'} = {\bf K}^\hbar_k$ uses only the exchange matrix $\varepsilon$ (and $k$), rather than the full topological information  of $\Delta$. Notice that the exchange matrix $\varepsilon=\varepsilon_\Delta$ of a labeled ideal triangulation $\Delta$ is insensitive to the mapping class group action, namely
$$
\varepsilon_\Delta = \varepsilon_{\h.\Delta}.
$$
That is, for any $i,j \in I$, the $(i,j)$-th entry of $\varepsilon_\Delta$ coincides with the $(i,j)$-th entry of $\varepsilon_{\h.\Delta}$. From this we deduce the invariance
$$
{\bf K}^\hbar_{\Delta,\Delta'} = {\bf K}^\hbar_{\h.\Delta,\h.\Delta'}
$$
in the case when $(\Delta\leadsto \Delta')$ is a flip. In the case when $(\Delta\leadsto \Delta')$ is a label permutation, i.e. $\Delta' = P_\sigma(\Delta)$ for some permutation $\sigma$ of $I$, it is easy to see that this invariance again holds. It then follows that the invariance holds for an arbitrary change of labeled ideal triangulations $(\Delta\leadsto \Delta')$.

\vs

Consider the quotient groupoid
$$
{\rm PCG}(\frak{S}) := {\rm Pt}(\frak{S})/{\rm MCG}(\frak{S})
$$
and call it the {\it Ptolemy class groupoid}. An object of ${\rm PCG}(\frak{S})$ representing $\Delta \in {\rm Ob}({\rm Pt}(\frak{S}))$ is denoted by $[\Delta]$, and a morphism of ${\rm PCG}(\frak{S})$ representing a morphism $(\Delta\leadsto \Delta')$ of ${\rm Pt}(\frak{S})$ is denoted by $[\Delta\leadsto \Delta']$. In particular, we have
$$
[\Delta] = [\h.\Delta] \quad\mbox{and}\quad
[\Delta\leadsto \Delta'] = [\h.\Delta \leadsto \h.\Delta']
$$
for each $\h\in {\rm MCG}(\frak{S})$. Using this notation, one can define
$$
\mathscr{H}_{[\Delta]} := \mathscr{H}_\Delta \quad\mbox{and}\quad 
{\bf K}^\hbar_{[\Delta\leadsto\Delta']} := {\bf K}^\hbar_{\Delta,\Delta'},
$$
which are well-defined by the invariance statements we have above. Thus the results of \S\ref{sec:the_unitary_intertwiners} yield a contravariant functor
\begin{align}
\label{eq:PCG_functor}
{\rm PCG}(\frak{S}) \to {\rm Hilb}
\end{align}
which sends an object $[\Delta]$ to the Hilbert space $\mathscr{H}_{[\Delta]}$ and a morphism $[\Delta \leadsto \Delta']$ to ${\bf K}^\hbar_{[\Delta\leadsto \Delta']}$.

\subsection{A trilogy of representations of the mapping class group}
\label{subsec:trilogy_of_representations}

Pick any labeled ideal triangulation $\Delta$ of $\frak{S}$, and consider the corresponding object $[\Delta]$ of the Ptolemy class groupoid ${\rm PCG}(\frak{S}) = {\rm Pt}(\frak{S})/{\rm MCG}(\frak{S})$. Note that for each mapping class $\h\in {\rm MCG}(\frak{S})$, the morphism $(\Delta \leadsto \h.\Delta)$ of the Ptolemy groupoid ${\rm Pt}(\frak{S})$ yields a morphism $[\Delta\leadsto \h.\Delta]$ of ${\rm PCG}(\frak{S})$ from $[\Delta]$ to $[\h.\Delta]=[\Delta]$. Hence $[\Delta\leadsto \h.\Delta]$ is an automorphism of $[\Delta]$. This gives a map
$$
\mathcal{F}_{\Delta} ~:~ {\rm MCG}(\frak{S}) \to {\rm Aut} ([\Delta]), \qquad  \h \mapsto [\Delta \leadsto \h.\Delta]
$$
For $\h_1,\h_2 \in {\rm MCG}(\frak{S})$ note
\begin{align*}
\mathcal{F}_\Delta(\h_1) \circ \mathcal{F}_\Delta(\h_2) 
& = [\Delta\leadsto \h_1.\Delta] \circ [\Delta\leadsto \h_2.\Delta] \\
& = [\h_2.\Delta \leadsto \h_2.(\h_1.\Delta)] \circ [\Delta\leadsto \h_2.\Delta] \\
& = [\Delta \leadsto \h_2.(\h_1.\Delta)] \\
& = [\Delta \leadsto (\h_2 \h_1).\Delta] \\
& = \mathcal{F}_\Delta(\h_2 \h_1),
\end{align*}
so $\mathcal{F}_\Delta$ is an anti-homomorphism of groups. It is easy to see that unless $\frak{S}$ is a once-punctured torus, $\mathcal{F}_\Delta$ is an isomorphism, since the action of ${\rm MCG}(\frak{S})$ on ${\rm Ob}({\rm Pt}(\frak{S}))$ is free. When $\frak{S}$ is a once-punctured torus, the kernel $\mathcal{F}_\Delta$ is generated by the `hyperelliptic' involution which is of order 2.

\vs

Composing $\mathcal{F}_\Delta$ with the functor in eq.\eqref{eq:PCG_functor}, one obtains a map
\begin{align*}
{\renewcommand{\arraystretch}{1.4} \begin{array}{rcl}
\rho^\hbar_\Delta ~:~ {\rm MCG}(\frak{S}) & \to &  {\rm U}(\mathscr{H}_{[\Delta]}) = \{\mbox{unitary maps $\mathscr{H}_{[\Delta]} \to \mathscr{H}_{[\Delta]}$}\}, \\
\h & \mapsto & {\bf K}^\hbar_{\mathcal{F}_\Delta(\h)} = {\bf K}^\hbar_{[\Delta\leadsto \h.\Delta]} = {\bf K}^\hbar_{\Delta,\h.\Delta} :  \mathscr{H}_{[\Delta]} \equiv \mathscr{H}_{\h.\Delta} \to \mathscr{H}_{\Delta} = \mathscr{H}_{[\Delta]}. 
\end{array}}
\end{align*}
It is easy to see that $\rho^\hbar_\Delta$ is a group homomorphism. Just to write down explicitly, note that for $\h_1,\h_2\in {\rm MCG}(\frak{S})$, one has
\begin{align*}
\rho^\hbar_\Delta(\h_1) \circ \rho^\hbar_\Delta(\h_2) 
& = {\bf K}^\hbar_{[\Delta\leadsto \h_1.\Delta]} \circ {\bf K}^\hbar_{[\Delta \leadsto \h_2.\Delta]} \\
& = {\bf K}^\hbar_{[\Delta\leadsto \h_1.\Delta]} \circ {\bf K}^\hbar_{[\h_1.\Delta \leadsto \h_1.(\h_2.\Delta)]} \\
& = {\bf K}^\hbar_{\Delta,\h_1.\Delta} \circ {\bf K}^\hbar_{\h_1.\Delta, (\h_1 \h_2).\Delta} \\
& = {\bf K}^\hbar_{\Delta,(\h_1 \h_2).\Delta} \qquad \qquad\qquad (\mbox{$\because$ eq.\eqref{eq:K_consistency}}) \\
& = {\bf K}^\hbar_{[\Delta\leadsto (\h_1 \h_2).\Delta]} \\
& = \rho^\hbar_\Delta(\h_1 \h_2).
\end{align*}
This $\rho^\hbar_\Delta$ is the promised unitary representation of ${\rm MCG}(\frak{S})$ on the Hilbert space $\mathscr{H}_{[\Delta]} = L^2(\mathbb{R}^I)$ coming from the quantization of the moduli space $\mathcal{GH}_\Lambda(\frak{S}\times \mathbb{R})$ of three-dimensional gravity. Roughly summarizing, for each mapping class $\h \in {\rm MCG}(\frak{S})$ and a labeled ideal triangulation $\Delta$, one finds a sequence of flips $\mu_k$ and label permutations $P_\sigma$ connecting $\Delta$ and $\h.\Delta$, and takes the composition of the intertwiners ${\bf K}^\hbar_k$ and ${\bf P}_\sigma$ corresponding to these flips and label permutations, to construct a unitary operator $\rho^\hbar_\Delta(\h)$. Recall that the flip intertwiner ${\bf K}^\hbar_k$ is a composition of a simple operator on $L^2(\mathbb{R}^I)$ induced by a linear automorphism of $\mathbb{R}^I$ and an operator expressed by the three kinds of quantum dilogarithm functions, depending on the cosmological constant $\Lambda$. 

\vs

Let us briefly discuss the dependence of this representation $\rho^\hbar_\Delta$ on the choice of a labeled ideal triangulation $\Delta$. To any ordered pair of labeled ideal triangulations $\Delta$ and $\Delta'$, is associated a unique unitary intertwiner ${\bf K}^\hbar_{\Delta,\Delta'} : \mathscr{H}_{\Delta'} \to \mathscr{H}_\Delta$, which can be interpreted as the unitary map ${\bf K}^\hbar_{[\Delta\leadsto \Delta']} : \mathscr{H}_{[\Delta']} \to \mathscr{H}_{[\Delta]}$. It is easy to observe that this map ${\bf K}^\hbar_{[\Delta\leadsto \Delta']}$ gives an intertwiner between the representations $\rho^\hbar_{\Delta'}$ and $\rho^\hbar_\Delta$, in the sense that the diagram
$$
\xymatrix{
\mathscr{H}_{[\Delta']} \ar[r]^-{\rho^\hbar_{\Delta'}(\h)} \ar[d]_{{\bf K}^\hbar_{[\Delta\leadsto \Delta']}} & \mathscr{H}_{[\Delta']} \ar[d]^{{\bf K}^\hbar_{[\Delta\leadsto \Delta']}} \\
\mathscr{H}_{[\Delta]} \ar[r]^-{\rho^\hbar_\Delta(\h)} & \mathscr{H}_{[\Delta]}
}
$$
commutes for all $\h\in {\rm MCG}(\frak{S})$, i.e. the following equation holds:
$$
{\bf K}^\hbar_{[\Delta\leadsto \Delta']} \circ \rho^\hbar_{\Delta'}(\h) = \rho^\hbar_\Delta(\h) \circ {\bf K}^\hbar_{[\Delta\leadsto \Delta']}.
$$
Hence, one can conclude that the representations $\rho^\hbar_\Delta$ on $\mathscr{H}_{[\Delta]}$ for different labeled ideal triangulations $\Delta$ are isomorphic to each other through a consistent system of isomorphisms ${\bf K}^\hbar_{[\Delta\leadsto \Delta']}$ between them.

\subsection{Irreducible representations satisfying quantum constraints}
\label{subsec:irreducible_representations}

The representation of ${\rm MCG}(\frak{S})$ constructed above is built on the unitary intertwiners ${\bf K}^\hbar_{\Delta,\Delta'}$, for which we used a particular self-adjoint operators ${\bf x}_i$ and ${\bf y}_i$ in eq.\eqref{eq:x_i_y_i_operators} on the Hilbert space $\mathscr{H}_\Delta$ for each $\Delta$. They were taken as forming a standard solution to the Heisenberg relations in eq.\eqref{eq:x_i_y_i_Heisenberg}. But they did not satisfy the quantum constraint equations in eq.\eqref{eq:quantum_constraint_equations}, so they cannot be regarded as honest quantizations of the coordinate functions $x_i$ and $y_i$ on the space $\mathcal{GH}_\Lambda(\frak{S} \times \mathbb{R})$. Another related fact is that they do not form an irreducible representation. In the present subsection we give a recipe to obtain an irreducible system of self-adjoint operators ${\bf x}_i$ and ${\bf y}_i$ (for $i\in \Delta$) on some Hilbert space satisfying both eq.\eqref{eq:x_i_y_i_Heisenberg} and eq.\eqref{eq:quantum_constraint_equations}. We note that this recipe is only hinted in \cite{KS}, but details are not given there.

\vs

The problem of finding such operators can be viewed as a problem of constructing representations of a generalized Heisenberg algebra, which was thoroughly studied in the literature; for example, one can find a version useful for us stated in \cite{Lion}. The construction of such representations as necessary for our purposes is reviewed in \cite{K21b}, which we apply to our situation here. Consider the vector space $(\mathbb{R}^\Delta)^2 = (\mathbb{R}^I)^2 = \mathbb{R}^{12g-12+6n}$ equipped with the skew-symmetric bilinear form $\langle \cdot, \cdot \rangle$ given by $\smallmattwo{0}{\varepsilon}{-\varepsilon}{0}$ on the basis $\{v_i,w_i\}_{i\in \Delta}$, i.e. $\langle v_i, v_j \rangle = 0 = \langle w_i,w_j\rangle$, $\langle v_i, w_j \rangle = \varepsilon_{ij} = \varepsilon_{i,j}$. Choose any Lagrangian subspace $\mathscr{L}$ of this space $((\mathbb{R}^\Delta)^2,\langle \cdot, \cdot \rangle)$, i.e. a maximal isotropic subspace; here, isotropic means $\langle v,w\rangle=0$, $\forall v,w\in \mathscr{L}$. Since $\dim (\ker \varepsilon) = n$ as mentioned in \S\ref{subsec:brief_overview}, one can observe that $\mathscr{L}$ is of dimension $2n+ \frac{1}{2}(12g-12+4n)$. Upon a choice of $\mathscr{L}$ and a basis $\mathcal{B}$ of the complement of $\mathscr{L}$ in $(\mathbb{R}^\Delta)^2$, using a basic harmonic analysis on the Heisenberg group for the Heisenberg Lie algebra defined as a $1$-dimensional central extension of the abelian Lie algebra $(\mathbb{R}^\Delta)^2$ constructed by using the form $\langle \cdot, \cdot \rangle$, an irreducible system of self-adjoint operators ${\bf x}_i$ and ${\bf y}_i$ (for $i\in \Delta$) satisfying eq.\eqref{eq:x_i_y_i_Heisenberg} and eq.\eqref{eq:quantum_constraint_equations} can be constructed on a Hilbert space isomorphic to $L^2(\mathbb{R}^\Delta/\mathscr{L}) = L^2(\mathbb{R}^{\mathcal{B}}) = L^2(\mathbb{R}^{6g-6+2n})$ with respect to the usual Lebesgue measure on $\mathbb{R}^{6g-6+2n}$. The resulting operators ${\bf x}_i$ and ${\bf y}_i$ are $\mathbb{R}$-linear combinations of the position and the momentum operators on $L^2(\mathbb{R}^{6g-6+2n}, \bigwedge_{j=1}^{6g-6+2n} ds_j)$. %Hence, the problem boils down to finding real numbers $a_j$, $b_j$, $c_j$ and $d_j$ for $j=1,\ldots,6g-6+2n$ such that the self-adjoint operators
%$$
%{\bf x}_i = \sum_{j=1}^{6g-6+2n} (a_j s_j + b_j {\rm i} \frac{\partial}{\partial s_j})
%\quad\mbox{and}\quad
%{\bf y}_i = \sum_{j=1}^{6g-6+2n} (c_j s_j + d_j {\rm i} \frac{\partial}{\partial s_j})
%$$
%satisfy eq.\eqref{eq:x_i_y_i_Heisenberg} and eq.\eqref{eq:quantum_constraint_equations}; this is just a simple linear algebra problem. 

\vs

Using the last comment, one can just try to look, by hand, for a system of $\mathbb{R}$-linear combinations of the position and the momentum operators that satisfy both eq.\eqref{eq:x_i_y_i_Heisenberg} and eq.\eqref{eq:quantum_constraint_equations}. But we could also use the full power of the construction we discussed above (reviewed in \cite{K21b}), which gives an explicit solution. For this, we first recall from \S\ref{subsec:brief_overview} that the (column) vectors $(v_{i,p})_{i\in \Delta}$ for punctures $p$ form a basis of the kernel of $\varepsilon$; in particular, $\sum_{j\in \Delta} \varepsilon_{ij} v_{j,p}=0$ for each puncture $p$ and each $i \in \Delta$. This means that $\{ \sum_{i \in \Delta} v_{i,p}v_i, \sum_{i \in \Delta} v_{i,p} w_i \}_{\mbox{\small $p$:punctures}}$ forms a basis of the kernel (or, radical) of $\langle \cdot, \cdot \rangle$. Since the $(6g-6+3n)\times n$ matrix $(v_{i,p})_{\mbox{\small $i\in \Delta$, $p$:punctures}}$ is of full rank $n$, one can find $n$ pivot rows, say row numbers $i_1,i_2,\ldots,i_n \in \Delta$. Then one can  transform this matrix into a column reduced echelon form using column operations, with leading 1 at the pivot rows, i.e. at $i_j$-th rows for $j=1,\ldots,n$. One can further assume that at each $i_j$-th row, there is only one non-zero entry, say at $p_j$-th place where the entry is $1$. Here $p_1,\ldots,p_n$ is an enumeration of (or, labels of) punctures. Denote this reduced matrix by $(w_{i,p})_{\mbox{\small $i\in \Delta$, $p$:punctures}}$; by the condition we imposed, for each $j=1,\ldots,n$, we have $w_{i_j,p_j}=1$, and $w_{i_j,p}=0$ for each $p\neq p_j$. Hence, if we write
$$
\mathring{\Delta} := \Delta \setminus\{i_1,\ldots,i_n\},
$$
which is a collection of $6g-6+2n$ arcs, then for each $j=1,\ldots,n$, in the $p_j$-th column we have $w_{i_{j'},p_j}=0$ whenever $j' \neq j$. Each column of the reduced matrix $(w_{i,p})_{i,p}$ is still in the kernel of the matrix $\varepsilon$, so $\sum_{k \in \Delta} \varepsilon_{ik} w_{k,p_j}=0$, and hence
\begin{align}
\label{eq:pivoted_kernel_basis}
\varepsilon_{i,i_j} + \textstyle \sum_{k \in \mathring{\Delta}} \varepsilon_{ik} w_{k,p_j}=0
\end{align}
holds for each $j=1,\ldots,n$ and $i\in \Delta$. The linear relations among the rows of $(v_{i,p})_{i,p}$ are the same as those of $(w_{i,p})_{i,p}$. For the latter matrix, for each $j=1,\ldots,n$, the $i_j$-th row is $p_j$-th standard basis row vector; hence, for each $i\in \mathring{\Delta}$, the $i$-th row $(w_{i,p})_{\mbox{\small $p$:punctures}}$ is the sum of $w_{i,p_j}$ times $i_j$-th row, where $j$ runs in $\{1,\ldots,n\}$. The same relation among the rows hold for the former matrix $(v_{i,p})$, too; for each $i \in \mathring{\Delta}$, one has
\begin{align}
\label{eq:reduced_matrix_row_relation}
v_{i,p} = \textstyle \sum_{j=1}^n w_{i,p_j} v_{i_j,p}
\end{align}
for each puncture $p$.

\vs

Now, for an example of the construction of Heisenberg algebra representation, let's choose the Lagrangian $\mathscr{L}$ of $(\mathbb{R}^\Delta)^2$ to be spanned by $w_i$ for $i\in \Delta$, together with $\sum_{i\in \Delta} v_{i,p} v_i$ for punctures $p$. As a basis $\mathcal{B}$ of the complement of $\mathscr{L}$ we take $v_i$ for $i \in \mathring{\Delta}$. Then, the construction reviewed in \cite{K21b} (especially, Examples 2.9 and 2.10 there are useful), with a slight scaling applied at the end for our purposes, yields the following self-adjoint operators 
\begin{align}
\label{eq:irreducible_solution}
\left\{ {\renewcommand{\arraystretch}{1.8} \begin{array}{ll}
\displaystyle {\bf x}_i = - \pi {\rm i} \frac{\partial}{\partial s_i}, & \forall i \in \mathring{\Delta} = \Delta\setminus\{i_1,\ldots,i_n\}, \\
\displaystyle {\bf x}_{i_j} = \underset{k\in \mathring{\Delta}}{\sum} \pi {\rm i} w_{k, p_j} \frac{\partial}{\partial s_k}, & \forall j=1,\ldots,n, \\
\displaystyle {\bf y}_i = \underset{j \in \mathring{\Delta}}{\sum} \varepsilon_{ij} s_j, &  \forall i \in \Delta,
\end{array}}
\right.
%\quad 
{\renewcommand{\arraystretch}{1.4} \begin{array}{l}
\mbox{on the Hilbert space} \\
\mbox{$\mathscr{H}_\Delta=L^2(\mathbb{R}^{\mathring{\Delta}}, \bigwedge_{i \in \mathring{\Delta}} ds_i)$} \\
\hspace{7,3mm} \mbox{$= L^2(\mathbb{R}^{6g-6+2n})$}
\end{array}}
\end{align}
One can now forget about the construction reviewed in \cite{K21b} and just take the above final result. One can check that they indeed satisfy the desired equations, as follows. It is easy to verify $[{\bf x}_i, {\bf x}_j]=0=[{\bf y}_i,{\bf y}_j]$ for all $i,j\in \Delta$, as well as $[{\bf x}_i, {\bf y}_j] = \pi {\rm i} \varepsilon_{ij}$ for all $i \in \mathring{\Delta}$ and $j \in \Delta$. For $j=1,\ldots,n$ and $k\in \Delta$, note that 
\begin{align*}
[{\bf x}_{i_j}, {\bf y}_k] & =\textstyle [\sum_{a \in \mathring{\Delta}} \pi {\rm i} w_{a,p_j} \frac{\partial}{\partial s_a}, \, \sum_{b \in \mathring{\Delta}} \varepsilon_{kb} s_b] \\
& = \textstyle \sum_{a\in \mathring{\Delta}} \sum_{b\in \mathring{\Delta}} \pi {\rm i} w_{a,p_j} \varepsilon_{kb} \delta_{a,b} \\
& = \textstyle \sum_{a\in \mathring{\Delta}} \pi {\rm i} w_{a,p_j} \varepsilon_{ka} \\
& = \pi {\rm i} \varepsilon_{i_j,k},
\end{align*}
where the last equality holds by eq.\eqref{eq:pivoted_kernel_basis}. So the Heisenberg relations in eq.\eqref{eq:x_i_y_i_Heisenberg} are satisfied. Moreover, for each puncture $p$, note that
\begin{align*}
\textstyle \sum_{i \in \Delta} v_{i,p} {\bf y}_i = \sum_{i \in \Delta}  v_{i,p} \sum_{j \in \mathring{\Delta}} \varepsilon_{ij} s_j
= \sum_{j \in \mathring{\Delta}} (\sum_{i \in \Delta} v_{i,p}  \varepsilon_{ij}) s_j = 0
\end{align*}
holds by eq.\eqref{eq:v_in_the_kernel}, and that
\begin{align*}
\textstyle \sum_{i\in \Delta} v_{i,p} {\bf x}_i & \textstyle = \sum_{i \in \mathring{\Delta}} v_{i,p} (-\pi {\rm i} \frac{\partial}{\partial s_i}) + \sum_{j=1}^n v_{i_j,p} \sum_{k \in \mathring{\Delta}} (\pi {\rm i} w_{k,p_j} \frac{\partial}{\partial s_k}) \\
& = \textstyle \sum_{i \in \mathring{\Delta}} (-\pi {\rm i}) (v_{i,p} - \sum_{j=1}^n v_{i_j,p} w_{i,p_j}) \frac{\partial}{\partial s_i} \\
& = 0,
\end{align*}
where the last equality holds by eq.\eqref{eq:reduced_matrix_row_relation}. Hence the quantum constraint equations in eq.\eqref{eq:quantum_constraint_equations} are also satisfied. One can also verify that these operators forms an irreducible system, and that any irreducible system of self-adjoint operators satisfying eq.\eqref{eq:x_i_y_i_Heisenberg} and eq.\eqref{eq:quantum_constraint_equations} is unitarily equivalent to the above one, using (generalized) Stone-von Neumann theorem; see e.g. \cite[\S14]{Hall}. 

\vs

The choice of a pair of Lagrangian subspace $\mathscr{L}$ and a complementary basis $\mathcal{B}$ made above is not being claimed to be canonical. A different choice yields in general a different set of operators on a Hilbert space isomorphic to $L^2(\mathbb{R}^{6g-6+2n})$. The construction reviewed in \cite{K21b} in fact also gives an explicit unitary equivalence between different solutions, using basic harmonic analysis. Due to the irreducibility, the unitary equivalence is unique up to scalar, and is what is referred to as the {\it (Segal-)Shale-Weil intertwiner} \cite{Shale, Segal, Weil, Lion}, which can be thought of as a generalization of the Fourier transformation. The monomial transformation part operator ${\bf K}'_k$ we described in the present article can also be thought of as a special example of a Shale-Weil intertwiner for a certain setting.

\vs

If one decides to use the operators ${\bf x}_i$ and ${\bf y}_i$ as in eq.\eqref{eq:irreducible_solution} instead of eq.\eqref{eq:x_i_y_i_operators}, one should also modify the subsequent constructions in \S\ref{sec:the_unitary_intertwiners}--\S\ref{sec:MCG_representations} accordingly, in particular the unitary intertwiners ${\bf K}^\hbar_{\Delta,\Delta'}$. The automorphism part operator ${\bf K}^{\hbar \sharp}_k$ for a flip $\mu_k$ is given by the same formula in eq.\eqref{eq:K_sharp_formula}, this time using new ${\bf x}_k$ and ${\bf y}_k$. The monomial transformation part operator ${\bf K}'_k$ for a flip $\mu_k$, as well as the index-permutation operator ${\bf P}_\sigma$, must be constructed afresh. Due to the irreducibility of the new set of operators ${\bf x}_i$ and ${\bf y}_i$, these sought-for operators ${\bf K}'_k$ and ${\bf P}_\sigma$ are completely determined, up to scalar, by their conjugation actions as follows (and these equations for conjugation actions are what we really need from the intertwiners):
\begin{align*}
& \hspace{-3mm} {\bf K}'_k \, {\bf x}_i'\, ({\bf K}'_k)^{-1} = \left\{
\begin{array}{ll}
-{\bf x}_k & \mbox{if $i=k$}, \\
{\bf x}_i + [\varepsilon_{ik}]_+ {\bf x}_k & \mbox{if $i \neq k$},
\end{array}
\right.
\quad
{\bf K}'_k \, {\bf y}_i'\, ({\bf K}'_k)^{-1} = \left\{
\begin{array}{ll}
-{\bf y}_k & \mbox{if $i=k$}, \\
{\bf y}_i + [\varepsilon_{ik}]_+ {\bf y}_k & \mbox{if $i \neq k$},
\end{array}
\right. \\
& \hspace{-3mm} {\bf P}_\sigma {\bf x}'_{\sigma(i)} {\bf P}_\sigma^{-1} = {\bf x}_i, \hspace{44mm}
{\bf P}_\sigma {\bf y}'_{\sigma(i)} {\bf P}_\sigma^{-1} = {\bf y}_i, \quad \forall i \in I.
\end{align*}
In order to explicitly construct such operators, one could again use the above construction reviewed in \cite{K21b}, this time for the Shale-Weil intertwiners. To say briefly, take the case of ${\bf K}'_k : \mathscr{H}_{\Delta'} \to \mathscr{H}_\Delta$, where $\Delta' = \mu_k(\Delta)$. The construction of a Heisenberg representation for $\Delta$ has required us to choose some data $(\mathscr{L},\mathcal{B})$ of a Lagrangian subspace $\mathscr{L}$ of $((\mathbb{R}^\Delta)^2,\langle \cdot,\cdot\rangle)$ and a complementary basis $\mathcal{B}$. Likewise, that for $\Delta'$  needed a choice $(\mathscr{L}',\mathcal{B}')$ for the space $( (\mathbb{R}^{\Delta'})^2, \langle \cdot, \cdot \rangle')$. One then takes a pullback of $(\mathscr{L}',\mathcal{B}')$ along the map $C: (\mathbb{R}^{\Delta})^2 \to (\mathbb{R}^{\Delta'})^2$ given on the basis by $C(v_k) = -v_k'$, $C(w_k)=-w_k'$, $C(v_i) = v_i' + [\varepsilon_{ik}]_+ v_k'$ and $C(w_i) = w_i' + [\varepsilon_{ik}]_+ w_k'$ for $i\neq k$. This pullback $(C^{-1}(\mathscr{L}'), C^{-1}(\mathcal{B}'))$ forms a data of a Lagrangian and a complementary basis for $(\mathbb{R}^\Delta)^2$, hence yields a Heisenberg algebra representation. Then ${\bf K}'_k$ is essentially defined as the Shale-Weil intertwiner between this representation for $(C^{-1}(\mathscr{L}'), C^{-1}(\mathcal{B}'))$ and that for $(\mathscr{L},\mathcal{B})$.  For more details about how to do this, we refer the readers to \cite{K21b}, where such a task was carefully done for quantization of the Teichm\"uller space. It will be a good exercise for a reader to explicitly construct ${\bf K}'_k$ and ${\bf P}_\sigma$ for our setting, mimicking \cite{K21b}.

\vs

To conclude, the quantization result deduced from the (irreducible) system of operators satisfying eq.\eqref{eq:x_i_y_i_Heisenberg} and eq.\eqref{eq:quantum_constraint_equations}, as appearing in eq.\eqref{eq:irreducible_solution} for example which we discussed in the present subsection, instead of the reducible solution in eq.\eqref{eq:x_i_y_i_operators} which does not satisfy eq.\eqref{eq:quantum_constraint_equations}, is what could really be regarded as a sought-for honest quantization of the moduli space $\mathcal{GH}_\Lambda(\frak{S}\times \mathbb{R})$ for 3d gravity. 

\vs

One worthwhile comment is that the intertwiners ${\bf K}^\hbar_{\Delta,\Delta'}$ constructed in the present subsection for the irreducible representations may satisfy the consistency equations in eq.\eqref{eq:K_consistency} only up to multiplicative constants. This means that the corresponding ${\rm MCG}(\frak{S})$ representations $\rho^\hbar_\Delta$ constructed using the methods of 
\S\ref{subsec:Ptolemy_groupoid}--\S\ref{subsec:trilogy_of_representations} a priori would only be {\it projective} representations, i.e. $\rho^\hbar_\Delta(\h_1) \rho^\hbar_\Delta(\h_2) = \rho^\hbar_\Delta(\h_1\h_2)$ holds only up to multiplicative constants. See \cite{K21b} for a method of explicit computation of these constants as well as for their meaning.

\subsection{Meaning of the results}

As mentioned in \cite{KS}, it has essentially been shown in the literature that the moduli space $\mathcal{GH}_\Lambda(\frak{S}\times \mathbb{R})$ of 3d spacetimes, with the gravitational symplectic structure, is symplectomorphic to the cotangent bundle $T^*\mathcal{T}(\frak{S})$ with the canonical (contangent-bundle) symplectic structure, where $\mathcal{T}(\frak{S})$ is the {\it (cusped) Teichm\"uller space} of $\frak{S}$, which is the set of all finite-area complete hyperbolic metrics on $\frak{S}$, considered up to pullback by self-diffeomorphisms of $\frak{S}$ isotopic to identity. This is true regardless of the value of the cosmological constant $\Lambda \in \{-1,0,1\}$. 

\vs

So, one could say that the problem of quantization of $\mathcal{GH}_\Lambda(\frak{S}\times \mathbb{R})$ boils down to quantization of the cotangent bundle $T^*\mathcal{T}(\frak{S})$. As quantization of the cotangent bundle of a manifold has been studied extensively, one could just think of applying the known results. For example, once we choose a coordinate system on the base manifold, which is $\mathcal{T}(\frak{S})$ in our case, there is a standard canonical quantization of the coordinate functions on the base manifold and the corresponding coordinates of the cotangent fiber, as self-adjoint operators on the Hilbert space $L^2(\mathcal{T}(\frak{S}))$; these operators are just the position and the momentum operators. Even if one takes this well-known construction as a partial solution to the quantization problem, one has to do more. Namely, what are other classical observable functions on $T^*\mathcal{T}(\frak{S})$ that one should and could quantize, and what operators should be assigned to them, in a canonical manner?

\vs

The quantization constructed in \cite{KS} takes a completely different route. It does not make use of the identification of $\mathcal{GH}_\Lambda(\frak{S}\times \mathbb{R})$ with the cotangent bundle $T^*\mathcal{T}(\frak{S})$ and its canonical quantization (of its Darboux coordinate functions). Instead, based on the special coordinate systems on $\mathcal{GH}_\Lambda(\frak{S}\times \mathbb{R})$ valued in the ring
$
\mathbb{R}_\Lambda = \mathbb{R}[\ell]/(\ell^2=-\Lambda)
$
studied by Meusburger and Scarinci \cite{MS} using 3d geometry, a certain set of classical observables to be quantized is suggested, mimicking the upper cluster algebra from the theory of cluster algebras \cite{FZ} and cluster varieties \cite{FG06, FG09}. An algebraic deformation quantization of these observables is constructed, by translating the results from quantization of cluster varieties from surfaces \cite{AK} with suitable modifications. A corresponding operator representation on a Hilbert space is constructed, by developing an $\mathbb{R}_\Lambda$-version of Fock and Goncharov's representation of quantum cluster variety at $\mathbb{R}_{>0}$ \cite{FG09}. By `composing' the algebraic deformation quantization and the operator representation, one obtains a quantization of the special classical observables we chose by operators on a Hilbert space.

\vs

What is reviewed in the present article is the Hilbert space representation (without the algebraic deformation quantization part). Although this representation is described in \cite{KS} using a version of functional calculus valued in $\mathbb{R}_\Delta$, here we avoided working with the $\mathbb{R}_\Lambda$-version, and translated the final results without using $\mathbb{R}_\Lambda$. An advantage of avoiding $\mathbb{R}_\Lambda$ is that the results can be described using usual theory of self-adjoint operators and their functional calculus. A slight disadvantage is that the answers for different values of the cosmological constant $\Lambda$ look different from each other, not a priori directly related. If one is willing to work with the $\mathbb{R}_\Lambda$-version as in \cite{KS}, then one can understand the quantization results in a uniform manner for all values of $\Lambda$, maybe with some parts remaining only heuristic. Especially, the formula of eq.\eqref{eq:K_sharp_formula} for the automorphism part operator ${\bf K}^{\hbar\sharp}_k$ for the unitary intertwiner ${\bf K}^\hbar_k$ for a flip $\mu_k$ can be described by a single formula
$$
\Phi^{\ell\hbar}({\bf x}_k + \hat{\ell} \hbar {\bf y}_k) \, \Phi^{-\ell\hbar}({\bf x}_k + \hat{\ell} \hbar {\bf y}_k)
$$
using the somewhat heuristic $\mathbb{R}_\Lambda$-version $\Phi^{\ell\hbar}$ of the quantum dilogarithm, where $\hat{\ell}$ is a suitable operator that quantizes the element $\ell$ of $\mathbb{R}_\ell$; in particular, $\hat{\ell}^2 = - \Lambda \cdot {\rm id}$. We refer the readers to \cite{KS} for more details.

\vs

In any case, the quantization of $\mathcal{GH}_\Lambda(\frak{S}\times \mathbb{R})$ obtained in \cite{KS}, which we reviewed in the present article, is not based on the canonical quantization of the cotangent bundle $T^*\mathcal{T}(\frak{S})$. It uses 3d geometry and does depend on the value of $\Lambda$. As an important consequence, it also yields a trilogy of unitary representations of the mapping class group ${\rm MCG}(\frak{S})$ for the three values of $\Lambda$, as we are focusing in the present article. As a matter of fact, the natural classical actions of ${\rm MCG}(\frak{S})$ on $\mathcal{GH}_\Lambda(\frak{S}\times \mathbb{R})$ for different $\Lambda$ yield different actions on the cotangent bundle $T^*\mathcal{T}(\frak{S})$. So, in fact, our quantizations for different $\Lambda$ can be interpreted as giving quantizations of three different quotients of $T^*\mathcal{T}(\frak{S})$ by ${\rm MCG}(\frak{S})$.

\vs

When $\Lambda \neq 0$, the quantum dilogarithm function used in the quantization, $\Phi^{\pm \hbar}$ or $\Phi^{\pm {\rm i} \hbar}$, has the real parameter $\hbar$ in it, so the resulting ${\rm MCG}(\frak{S})$-representation $\rho^\hbar_\Delta$, based on the unitary intertwiners ${\bf K}^\hbar_{\Delta,\Delta'}$, can be thought of as forming a one-real-parameter family of representations. However, for the case $\Lambda=0$, there is no dependence on the quantum parameter $\hbar$, as observed in eq.\eqref{eq:K_sharp_formula}; so we get a single representation, not a family, although this is still a `quantum' representation.

\vs

As mentioned in \S\ref{subsec:flip_intertwiner}, for $\Lambda=-1$ which corresponds to the 3d anti-de Sitter geometry, the intertwiners and the resulting unitary ${\rm MCG}(\frak{S})$-representations obtained in \cite{KS} as written in the present article coincide with Fock and Goncharov's representation for the quantum cluster symplectic double variety of the Teichm\"uller space \cite{FG09}. The crucial technical ingredient is the non-compact quantum dilogarithm $\Phi^{\pm \hbar}$ of Faddeev and Kashaev \cite{Faddeev, FK94}, which actually goes back to Barnes 100 years ago \cite{B01}, and which is used already in the quantization of Teichm\"uller spaces by Kashaev \cite{Kash98} and by Chekhov and Fock \cite{CF, F} in 1990's. We note that the Chekhov-Fock quantization of Teichm\"uller spaces \cite{CF, F} is roughly `half' of the quantization described in the present article for $\Lambda=-1$, or equivalently of the quantization of \cite{FG09}; for Chekhov and Fock's intertwiner ${\bf K}^\hbar_k$ for a flip $\mu_k$, the automorphism part ${\bf K}^{\hbar\sharp}_k$ is given roughly by a single non-compact quantum dilogarithm $\Phi^\hbar({\bf x}_k + \hbar {\bf y}_k)$, which is a `half' of our answer for $\Lambda=-1$, and the monomial transformation part ${\bf K}'_k$ is given by a suitable Shale-Weil intertwiner, as explicitly written in \cite{K21b} for the irreducible representation.

\vs

For $\Lambda=1$ which corresponds to the 3d de Sitter geometry, the intertwiners and unitary ${\rm MCG}(\frak{S})$-representations are based on the version of a quantum dilogarithm function denoted by $\Phi^{\pm {\rm i}\hbar}$ \cite{KS}, which is a certain `modular double' version of the compact quantum dilogarithm $\psi^q(z) = \prod_{n \ge 0} (1+q^{2n+1}z)^{-1}$. We note that the compact quantum dilogarithm $\psi^q$ (of Faddeev and Kashaev \cite{Faddeev, FK94}) is not suitable for quantization of Teichm\"uller spaces because it does not make sense when $|q|=1$. This is the reason why the non-compact quantum dilogarithm $\Phi^{\pm \hbar}$ was used in quantum Teichm\"uller theory, which can be viewed as a {\it limit} of a certain ratio of two compact quantum dilogarithm functions, but not a genuine ratio of such. Meanwhile, our function $\Phi^{\pm {\rm i}\hbar}$, which does not seem to have been particularly considered in the literature before \cite{KS}, is a genuine ratio of two compact quantum dilogarithm functions, not a limit of such. A single copy of $\Phi^{\pm {\rm i}\hbar}$ does not provide a unitary operator, and only a suitable product of the two versions $\Phi^{{\rm i} \hbar}$ and $\Phi^{-{\rm i} \hbar}$ does, where these two copies represent the $\mathbb{R}_\Lambda$-holomorphic and the $\mathbb{R}_\Lambda$-anti-holomorphic parts respectively. In the end, the intertwiner for a flip is a unitary intertwiner made up by four genuine compact quantum dilogarithm functions $\psi^{\bf q}$ with real quantum parameters ${\bf q}$ related to each other in a specific way. As far as the author knows, this was the first time that honest compact quantum dilogarithm functions (without taking their limits) were used directly to give unitary representations of ${\rm MCG}(\frak{S})$ on a Hilbert space.

\vs

The intertwiners and unitary ${\rm MCG}(\frak{S})$-representations for $\Lambda=0$ which corresponds to the 3d Minkowski (`flat') geometry can be thought of as being built on the classical logarithm function, not quantum nor dilogarithm. They can also be regarded as coming from a `derivative' of the family of quantum dilogarithm functions, due to the nature of the ring $\mathbb{R}_\Lambda$ when $\Lambda=0$; indeed, notice that for an `analytic' function $f$, observe that $f(x+\ell y) = f(x) + \ell f'(x)y$ holds in the ring $\mathbb{R}_0 = \mathbb{R}[\ell]/(\ell^2=0)$, for $x,y\in \mathbb{R}$.

\section{Intertwiners and quantum group representations}

Here we discuss a relationship with the quantum group representations. For this we will first review the results of the joint work of the author with Igor Frenkel \cite{FK}.

\subsection{Quantum Teichm\"uller theory and quantum plane algebra}

Recall that to each labeled ideal triangulation $\Delta$ is associated a Hilbert space $\mathscr{H}_\Delta$, and to each change of labeled ideal triangulations $(\Delta\leadsto \Delta')$ is associated a unitary intertwining operator ${\bf K}^\hbar_{\Delta,\Delta'}$. An elementary case is when $\Delta'$ is obtained from $\Delta$ by a flip $\mu_k$ at an arc $k$, and we denoted the unitary intertwiner as ${\bf K}^\hbar_k$. One property of this flip intertwiner is that it satisfies the pentagon identity, as in eq.\eqref{eq:pentagon_operator_identity} of Thm.\ref{thm:consistency_equations_of_intertwiners}; a slight nuisance is that this version of the pentagon identity involves not just the flip intertwiners, but also an index permutation operator ${\bf P}_\sigma$. Analogous results hold for the Chekhov-Fock quantization of Teichm\"uller spaces \cite{CF, F}, which is also based on labeled ideal triangulations. The unitary intertwiner for a flip in this quantization also satisfies the pentagon identity as in eq.\eqref{eq:pentagon_operator_identity}, up to a multiplicative scalar. Any reader who is interested in this scalar can consult \cite{K21b, FS}.

\vs

Meanwhile, Kashaev's quantization of Teichm\"uller spaces \cite{Kash98} is based not on labeled ideal triangulations, but on ideal triangulations together with a different kind of extra data. It uses a {\it dotted ideal triangulation}, which consists of an ideal triangulation $\Delta$, a labeling of ideal triangles of $\Delta$, and the choice of a distinguished corner for each ideal triangle of $\Delta$ which is depicted by a dot $\bullet$ in pictures. To each dotted ideal triangulation, is associated a Hilbert space $\underset{\mbox{\small $t$: triangles}}{\bigotimes} L^2(\mathbb{R}) \cong L^2(\mathbb{R}^{\{\mbox{\small triangles}\}})=L^2(\mathbb{R}^{4g-4+2n})$ (here, the symbol $\bigotimes$ denotes the Hilbert space tensor product), and to each change of dotted ideal triangulations is associated a unitary intertwiner. There are three kinds of elementary changes: 1) $T_{ts}$ which is a flip at an arc adjacent to triangles $t$ and $s$ with a special dot configuration on these triangles, 2) $A_t$ which is a counterclockwise change of a dot in the triangle $t$, and 3) $P_\sigma$ which is a label permutation of triangles. The corresponding unitary intertwiners can be denoted by ${\bf T}^\hbar_{ts}$, ${\bf A}^\hbar_t$ and ${\bf P}_\sigma$. One of the consistency equations satisfied by these operators is another version of a pentagon identity, which reads
$$
{\bf T}^\hbar_{ts} {\bf T}^\hbar_{rt} = {\bf T}^\hbar_{rt} {\bf T}^\hbar_{rs} {\bf T}^\hbar_{ts}.
$$
The subscripts indicate which tensor factors of $\underset{\mbox{\small triangles}}{\bigotimes} L^2(\mathbb{R})$ are involved; more precisely, there is a single operator ${\bf T}^\hbar : L^2(\mathbb{R}) \otimes L^2(\mathbb{R}) \to L^2(\mathbb{R}) \otimes L^2(\mathbb{R})$ such that ${\bf T}^\hbar_{ij}$ means that this operator ${\bf T}^\hbar$ is being applied to the $i$-th and the $j$-th tensor factors $L^2(\mathbb{R})$. The above pentagon identity involves only three tensor factors, and its prototypical version can be written as
\begin{align}
\label{eq:clean_pentagon_identity}
{\bf T}^\hbar_{23} {\bf T}^\hbar_{12} = {\bf T}^\hbar_{12} {\bf T}^\hbar_{13} {\bf T}^\hbar_{23}
\end{align}
which is on $L^2(\mathbb{R})\otimes L^2(\mathbb{R})\otimes L^2(\mathbb{R})$. We note that this version of the pentagon identity involves only ${\bf T}^\hbar$, but not other operators like index permutation operators.

\vs

The pentagon identity naturally arises in several areas of mathematics, although maybe in slightly different forms. For example, the version as in eq.\eqref{eq:clean_pentagon_identity} acting on a triple tensor product is observed in a tensor (or, monoidal) category as a consistency axiom for the associativity morphism $(A\otimes B) \otimes C \to A\otimes (B\otimes C)$. Thus, it is natural to ask whether there exists a tensor category whose associativity morphisms recover the above operator ${\bf T}^\hbar$ appearing in Kashaev's quantum Teichm\"uller theory. This question is answered affirmative in the author's joint work with Igor Frenkel \cite{FK}, where we consider the category of certain well-behaved representations of a rather simple Hopf $*$-algebra called the {\it quantum plane} $\mathcal{B}_q$, which is an associative algebra generated by $X^{\pm 1}$ and $Y$, modulo the relation 
$$
XY = q^2 YX,
$$
where the associated quantum parameter $q = e^{\pi {\rm i} \hbar}$ is a complex number of modulus $1$. This algebra can be viewed as a quantum observable algebra of the plane $\mathbb{R}^2 = \{(x,y) : x,y\in \mathbb{R}\}$ with a standard symplectic structure $dx\wedge dy$, where $X$ and $Y$ are quantum counterparts of the classical functions $e^x$ and $e^y$. What may distinguish this algebra $\mathcal{B}_q$ from what is called the quantum torus is the $*$-structure, which is given by $X^*=X$ and $Y^*=Y$ (for the quantum plane), rather than $X^* = X^{-1}$ and $Y^* = Y^{-1}$ (for the quantum torus); this means that when we consider a representation of $\mathcal{B}_q$ on a Hilbert space, we would like $X$ and $Y$ to be represented by self-adjoint operators, rather than unitary operators. 

\vs

An irreducible `well-behaved' representation $\pi$ of $\mathcal{B}_q$ on a Hilbert space $\mathscr{H}$ is uniquely determined up to equivalence \cite{Schpaper}; it is given by $\mathscr{H} = L^2(\mathbb{R})$, where $X$ and $Y$ are represented by real exponentials of the position and the momentum operators. The algebra $\mathcal{B}_q$ can be viewed as a Borel subalgebra of the split real quantum group $\mathcal{U}_q(\frak{sl}(2,\mathbb{R}))$, and in particular is a Hopf algebra whose coproduct is given by $\Delta(X) = X\otimes X$ and $\Delta(Y) = Y\otimes X + 1 \otimes Y$ (here, $\Delta$ is the coproduct, and has nothing to do with an ideal triangulation!). Using the coproduct, one defines action of $\mathcal{B}_q$ on the tensor product of representations. We show that the tensor square $\mathscr{H}\otimes \mathscr{H}$ decomposes into direct integral $\int_\mathbb{R}^\oplus \mathscr{H}$ of the unique irreducible representations $\mathscr{H}$; we show this in the form of an isomorphism ${\bf F}^\hbar: \mathscr{H} \otimes \mathscr{H} \overset{\sim}{\longrightarrow} M \otimes \mathscr{H}$, where $M = L^2(\mathbb{R})$ is the trivial $\mathcal{B}_q$-module called the multiplicity module. We realize the decomposition isomorphism ${\bf F}^\hbar$ as a map $L^2(\mathbb{R})\otimes L^2(\mathbb{R}) \to L^2(\mathbb{R})\otimes L^2(\mathbb{R})$, given by an integral transformation expressed using the non-compact quantum dilogarithm function $\Phi^\hbar$; we note that it is also possible to express this operator ${\bf F}^\hbar$ as a composition of a Shale-Weil type operator and an operator given by functional calculus applied to the function $\Phi^\hbar$, just as we became familiar with through the present article. 

\vs

The tensor cube of $\mathscr{H}$ has two decompositions, $(\mathscr{H} \otimes \mathscr{H}) \otimes \mathscr{H}$ and $\mathscr{H} \otimes (\mathscr{H} \otimes \mathscr{H})$. The coassociatity of the coproduct implies that the `identity' map $(\mathscr{H} \otimes \mathscr{H}) \otimes \mathscr{H} \to \mathscr{H} \otimes (\mathscr{H} \otimes \mathscr{H})$ is an isomorphism of representations. Meanwhile, the parentheses tell us how to decompose each of these representations into direct integral of $\mathscr{H}$, using the decomposition maps ${\bf F}^\hbar$. More precisely, we obtain the diagram
$$
\xymatrix{
(\mathscr{H} \otimes \mathscr{H}) \otimes \mathscr{H} \ar[r]^-{{\bf F}^\hbar_{12}} \ar[d]_{{\rm id}} & M\otimes \mathscr{H} \otimes \mathscr{H} \ar[r]^-{{\bf F}^\hbar_{23}} & M\otimes M \otimes \mathscr{H} \ar@{-->}[d] \\
\mathscr{H} \otimes (\mathscr{H} \otimes \mathscr{H})  \ar[r]^-{{\bf F}^\hbar_{23}} & \mathscr{H} \otimes M \otimes \mathscr{H} \ar[r]^-{{\bf F}^\hbar_{13}} & M\otimes M \otimes \mathscr{H}
}
$$
By composing the four ${\bf F}^\hbar$ in the diagram, one gets the rightmost vertical map ${\bf F}^\hbar_{13} {\bf F}^\hbar_{23} ({\bf F}^\hbar_{12})^{-1} ({\bf F}^\hbar_{23})^{-1} : M\otimes M \otimes \mathscr{H} \to M \otimes M \otimes \mathscr{H}$ which we prove to be of the form ${\bf T}^\hbar \otimes {\rm id}$ for some ${\bf T}^\hbar : M \otimes M \to M \otimes M$. So, the `identity' associativity morphism $(\mathscr{H} \otimes \mathscr{H}) \otimes \mathscr{H}  \to \mathscr{H} \otimes (\mathscr{H} \otimes \mathscr{H})$ is encoded via ${\bf F}^\hbar$ as a map ${\bf T}^\hbar : M\otimes M \to M \otimes M$ between the multiplicity modules. We note that, by using the pentagon operator identity satisfied by the non-compact quantum dilogarithm $\Phi^\hbar$, we computed ${\bf T}^\hbar$ to be again of a similar form as ${\bf F}^\hbar$, namely composition of a Shale-Weil type operator and an operator expressed by $\Phi^\hbar$.

\vs

Consider the five different possible decompositions of the quadruple tensor product $\mathscr{H}\otimes \mathscr{H} \otimes \mathscr{H} \otimes \mathscr{H}$, represented by the choices of paranthesizing such as $(((\mathscr{H} \otimes \mathscr{H}) \otimes \mathscr{H}) \otimes \mathscr{H}$. The identity maps between them can be translated into the maps ${\bf T}^\hbar$ between the multiplicity modules $M \otimes M \otimes M$. The identity maps all commute with each other, and this commutativity is translated into the pentagon identity of five copies of ${\bf T}^\hbar$, which is exactly of the form as in eq.\eqref{eq:clean_pentagon_identity}. 

\vs

Part of the main theorem of \cite{FK} says that this operator ${\bf T}^\hbar$ coming from the representation theory of the quantum plane algebra $\mathcal{B}_q$, or in fact more precisely of its modular double $\mathcal{B}_q \otimes \mathcal{B}_{q^\vee}$ with $q^\vee = e^{ \pi {\rm i} / \hbar}$, coincides with the unitary intertwiner ${\bf T}^\hbar$ for a flip $T_{ts}$ of dotted ideal triangulations appearing in Kashaev's quantum Teichm\"uller theory \cite{Kash98}, up to conjugation by a unitary map. Moreover, the rigid tensor category structure for the representations of $\mathcal{B}_q$ is also studied. By making use of the identification $M \cong {\rm Hom}_{\mathcal{B}_q}(\mathscr{H}, \mathscr{H} \otimes \mathscr{H})$ of $M$ with the space of intertwiners, and by permuting the roles of the three copies of $\mathscr{H}$ appearing in the right hand side of this isomorphism using various natural isomorphisms involving the dual representations, we construct a map ${\bf A}^\hbar : M \to M$ of order three, and we have shown that it essentially coincides with the intertwiner ${\bf A}^\hbar$ for a dot-change $A_t$ in Kashaev's quantum Teichm\"uller theory.

\subsection{3d quantum gravity and quantum pseudo-K\"ahler plane algebra}

Let's come back now to 3d quantum gravity. Inspired by the flip intertwiners ${\bf K}^\hbar_k$ for the quantization of $\mathcal{GH}_\Lambda(\frak{S}\times \mathbb{R})$ constructed in \cite{KS}, the author has developed results analogous to \cite{FK} suited to 3d quantum gravity, especially for the case of positive cosmological constant $\Lambda=1$, in \cite{KimQC}. First, a new Hopf $*$-algebra $\mathcal{C}_{\bf q}$ is introduced, where the associated quantum parameter ${\bf q} = e^{-\pi \hbar}$ is real; $\mathcal{C}_{\bf q}$ is generated by $Z_1^{\pm 1}$ and $Z_2$, with the relation
$$
Z_1 Z_2 = {\bf q}^2 Z_2 Z_1.
$$
This algebra structure seems just like the quantum torus algebra or the quantum plane algebra, this time with a real quantum parameter. What distinquishes $\mathcal{C}_{\bf q}$ vastly with these previous algebras is the $*$-structure. In fact, being more precise, $C_{\bf q}$ is generated by $Z_1^{\pm 1}$, $Z_2$, $(Z_1^*)^{\pm 1}$, $Z_2^*$, modulo the above relation, and also additional relations saying that each of $Z_1$ and $Z_2$ commutes with each of $Z_1^*$ and $Z_2^*$, and that $Z_1^* Z_2^* = {\bf q}^{-2} Z_2^* Z_1^*$. The notations of the generators immediately suggest the $*$-structure. The coproduct is given as before by $\Delta Z_1 = Z_1 \otimes Z_1$, $\Delta Z_2 = Z_2 \otimes Z_1 + 1 \otimes Z_2$, and similarly for $Z_1^*$ and $Z_2^*$. One can view this algebra as being a quantum observable algebra for the complex `plane' $\mathbb{C}^2 = \{(z_1,z_2) : z_1,z_2 \in \mathbb{C}\}$, equipped with the pseudo-K\"ahler metric $dz_1 \otimes d\ol{z}_2 + dz_2 \otimes d\ol{z}_1$. Hence this algebra is called the {\it quantum pseudo-K\"ahler plane} in \cite{KimQC}.

\vs

A basic representation of $\mathcal{C}_{\bf q}$, or more precisely of its modular double $\mathcal{C}_{\bf q} \otimes \mathcal{C}_{1/{\bf q}^\vee}$ (with ${\bf q}^{\vee} =  e^{-\pi/\hbar}$), is constructed on $\mathscr{H} = L^2(\mathbb{R}^2)$, where the generators are represented by normal operators given by exponentials of simple $\mathbb{C}$-linear combinations of the position and the momentum operators. Then all the constructions as in \cite{FK} are mimicked by studying the representation theory of $\mathcal{C}_{\bf q}$, yielding the new corresponding versions of the operators ${\bf F}^\hbar$, ${\bf T}^\hbar$ and ${\bf A}^\hbar$, which satisfy the same consistency relations as before, including the pentagon identity in eq.\eqref{eq:clean_pentagon_identity}. Notice that ${\bf F}^\hbar$ and ${\bf T}^\hbar$ this time are given in terms of the modular-double version of the compact quantum dilogarithm, i.e $\Phi^{\pm {\rm i}\hbar}$.

\vs

For the case of the quantum plane algebra $\mathcal{B}_q$ as in \cite{FK}, these operators ${\bf T}^\hbar$ and ${\bf A}^\hbar$ are what constitute Kashaev's version of quantization of Teichm\"uller spaces. Note that Kashaev's quantum Teichm\"uller theory \cite{Kash98} is based Kashaev's special coordinate functions which are ratios of Penner's lambda length coordinates on the `decorated' Teichm\"uller space $\mathcal{T}^d(\frak{S})$ \cite{P87}. Kashaev's functions do not really form a coordinate system of a Teichm\"uller space, and can be viewed as being a coordinate system of some related space, given by the product of Penner's projectivized decorated Teichm\"uller space $\mathcal{T}^d(\frak{S})/\mathbb{R}_{>0}$ and the first cohomology $H^1(\frak{S};\mathbb{R})$ \cite[Prop.5]{Kash98}. 

\vs

Similarly for the case of the quantum pseudo-K\"ahler plane algebra $\mathcal{C}_{\bf q}$ as in \cite{KimQC}, the operators ${\bf T}^\hbar$ and ${\bf A}^\hbar$ should be what constitute a Kashaev-type quantization of 3d gravity for $\Lambda=1$, which is never considered. One should define a version of `decorated' moduli space of 3d gravity, establish 3d complexified version of lambda lengths there (one may want to see \cite{NN}), consider certain ratios of them, and so on. Such a construction is still yet to be found and left open.

\vs

Note that in \cite{KimQC}, only the case $\Lambda=1$ is dealt with in detail, which inspired the new Hopf algebra as above. One may continue investigating the cases for the remaining two values of $\Lambda$. Most likely, the case of $\Lambda=-1$, which is treated in \cite[\S14]{KimQC} without proof, would just be a `symplectic double' version of what was obtained in \cite{FK} already. The case $\Lambda=0$ might be interesting. The first step would be to define a certain new version of `quantum plane' algebra, probably over $\mathbb{R}_0$ or over the suitable complexification of $\mathbb{R}_0$ (see \cite{KS}), which is without the quantum parameter at all! The representation theory of this algebra should yield the `flat' versions of the operators ${\bf T}$ and ${\bf A}$ as well, where ${\bf T}$ is given in terms of the usual logarithm function. Then again, the Kashaev-type 3d quantum gravity theory should also be sought for.

\section{Some applications and future directions}

\ul{Modular functor conjecture}

\vs

Notice that the quantum Teichm\"uller theory associates to each punctured surface $\frak{S}$ a family of unitary (projective) representations of ${\rm MCG}(\frak{S})$ on a Hilbert space $\mathscr{H}_\frak{S}$; the construction a priori depends on the choice of an ideal triangulation $\Delta$ of $\frak{S}$, but the results for different $\Delta$ are consistently identified with each other. The {\it (${\it SL}_2$) modular functor conjecture} speculates that this assignment of Hilbert spaces to different surfaces is organized in such a way that it forms a so-called 2-dimensional `modular functor' corresponding to some 2-dimensional conformal field theory. A crucial axiom to be checked is, for each punctured surface $\frak{S}$ and a simple loop $\gamma$ in $\frak{S}$ that does not contract to an interior point or a puncture, if $\frak{S}_\gamma$ denotes the new surface obtained from $\frak{S}$ by cutting along $\gamma$ (and shrinking the two holes thus created from $\gamma$ to punctures), then the Hilbert space representations $\mathscr{H}_\frak{S}$ and $\mathscr{H}_{\frak{S}_\gamma}$ associated to $\frak{S}$ and $\frak{S}_\gamma$ by quantum Teichm\"uller theory should be related through a special way. For details, see \cite{Teschner} and \cite{FG09}. A modified version is verified in \cite{Teschner}, based on a novel adaptation of  Kashaev's quantum Teichm\"uller. Still, a more direct version using the Chekhov-Fock-Goncharov quantum Teichm\"uller theory \cite{CF, F, FG09, K21b} has not been fully established in a rigorous manner, and some of the first steps toward it are performed in \cite{Kim24}. 

\vs

Meanwhile, the version of 3d quantum gravity constructed in \cite{KS} also associates a Hilbert space representation of ${\rm MCG}(\frak{S})$ to each punctured surface $\frak{S}$. Hence one can formulate and try to solve the corresponding `${\rm SL}_2(\mathbb{R}_\Lambda)$' version of the modular functor conjecture. In the end, it would yield a 2-dimensional modular functor corresponding to some 2-dimensional conformal field theory.

\vs

Meanwhile, notice that there is a general correspondence between 2-dimensional conformal field theories (2d CFT) and 3-dimensional topological quantum field theories (3d TQFT); see \cite{Funar}. So the above modular functor conjectures hint at certain 3d TQFT lurking behind the quantum Teichm\"uller theory and 3d quantum gravity theory. On the other hand, Kashaev and Andersen approach 3d TQFT more directly from quantum Teichm\"uller theory \cite{AndersenKashaev}. One could try an analog of this using the 3d quantum gravity theory of \cite{KS}.

\vspace{5mm}

\ul{The universal case}

\vs

The universal Teichm\"uller space has been known to be related to the universal phase space of 3d anti-de Sitter gravity; see \cite{BH, SK}. Thus, quantization of the universal Teichm\"uller space, which can be viewed as the Teichm\"uller space of the open unit disc, with certain restriction on the boundary behavior \cite{Penner93}, should be related to quantization of 3d anti-de Sitter gravity. Besides, the mapping class group for a version of the universal Teichm\"uller space is the `Ptolemy-Thompson group' $T$, which is one of the three infinite groups of R. Thompson \cite{CFP96}. In particular, quantization of universal Teichm\"uller space yields representations of $T$; see \cite{FS} and \cite{K16a}. On the other hand, Igor Frenkel and Bob Penner have suggested in \cite{FP} that the study of universal Teichm\"uller space might lead to the Monster 2d CFT which was constructed by I. Frenkel, J. Lepowsky and A. Meurman in 1980's \cite{FLM}. Keeping in mind Edward Witten's conjectural relationship between 3d quantum gravity and Monster 2d CFT \cite{Witten}, it should be worth investigating the universal case of the construction of \cite{KS}, and thinking how it could be connected to all these conjectures. We refer the readers to \cite{FK, K16a, FS, Penner93} for more details on how to approach a quantization of universal Teichm\"uller space.

\vspace{5mm}

\ul{Going to higher rank}

\vs

The Teichm\"uller space of a surface $\frak{S}$ can be thought of as some version (or, some component) of the representation variety ${\rm Hom}(\pi_1(\frak{S}), {\rm G})/{\rm G}$, when ${\rm G} = {\rm PSL}_2(\mathbb{R})$. The moduli space $\mathcal{GH}_\Lambda(\frak{S} \times \mathbb{R})$ of 3d gravity on $\frak{S} \times \mathbb{R}$ corresponds to the case when ${\rm G} = {\rm PSL}_2(\mathbb{R}_\Lambda)$ (see \cite{KS}). Meanwhile, the study of the representation variety for higher rank (real) Lie groups ${\rm G}$, as being a higher rank analog of the Teichm\"uller theory, is referred to as the {\it higher Teichm\"uller theory} \cite{FG06}. In some case like ${\rm G} = {\rm PGL}_3(\mathbb{R})$, the corresponding space can be interpreted as the moduli space of certain geometric structures on the surface $\frak{S}$, such as convex projective structures \cite{Goldman}. One might try to look for such a geometric description when we $\mathbb{R}_\Lambda$-complexify, e.g. in the case when ${\rm G} = {\rm PGL}_3(\mathbb{R}_\Lambda)$. Quantization of such $\mathbb{R}_\Lambda$-complexified moduli spaces, or their suitable cluster variety relatives, can in fact be obtained already using the results of \cite{KS}, because most of the constructions of \cite{KS} are formulated to work for general cluster varieties, not just for the space $\mathcal{GH}_\Lambda(\frak{S} \times \mathbb{R})$.

\vspace{5mm}

\ul{More general topology}

\vs

After the author's talk on \cite{KS} at Igor Frenkel's 70th birthday conference, Nikita Nekrasov in the audience asked what would happen to the case of 3d gravity when the topology of the surface changes (`in time'). The quantization methods of \cite{KS} depend on the $\mathbb{R}_\Lambda$-valued cluster coordinates developed in \cite{MS}, which work only on the space $\mathcal{GH}_\Lambda(\frak{S}\times \mathbb{R})$ of 3d gravity on the 3-manifold $\frak{S}\times \mathbb{R}$. For a moduli space of 3d gravity on a more general 3-manifold which does not look like a surface times an interval, the cluster variety structure has not been observed as far as the author knows, although somewhat similar looking coordinate change maps may have been found \cite{GTZ}. It seems that one might need to develop a new machinery to deal with these (e.g. to quantize them), instead of trying to apply the known theory of cluster varieties.

\section*{Acknowledgments}

This work was supported by the National Research Foundation of Korea (NRF) grant funded by the Korea government (MSIT) (No. 2020R1C1C1A01011151). H. Kim has been supported by a KIAS Individual Grant (MG047203, MG047204) at Korea Institute for Advanced Study. H. Kim thanks his advisor Igor B. Frenkel for giving continuous support and inspirations. H. Kim thanks the organizers of Igor's 70th birthday conference for providing a wonderful experience at the conference, and also for putting together this proceedings.

\bibliographystyle{amsalpha}

\begin{thebibliography}{FKV01}

\bibitem[AK17]{AK} D. G. L. Allegretti, H. Kim, {\it A duality map for quantum cluster varieties from surfaces}, Adv. Math. {\bf 306} (2017), 1164--1208.

\bibitem[AK14]{AndersenKashaev} J. E. Andersen, R. M. Kashaev, {\it A TQFT from quantum Teichm\"uller theory}, Commun. Math. Phys. {\bf 330} (2014), 887--934.


\bibitem[BK00]{BK} B. Bakalov, A. Kirillov Jr., {\it On the Lego-Teichm\"uller game}, Transform. Groups {\bf 5} (2000), 207--244.

\bibitem[B01]{B01} E. W. Barnes, {\it Theory of the double Gamma function}, Phil. Trans. Roy. Soc. A {\bf 196} (1901), 265--388.

\bibitem[BL07]{BL} F. Bonahon, X. Liu, {\it Representations of the quantum Teichm\"uller space and invariants of surface diffeomorphisms}, Geom. Topol. {\bf 11} (2007), 889--937.

\bibitem[BW11]{BW} F. Bonahon, H. Wong, {\it Quantum traces for representations of surface groups in ${\rm SL}_2(\mathbb{C})$}, Geom. Topol. {\bf 15} (2011), 1569--1615.

\bibitem[BH86]{BH} J. D. Brown, M. Henneaux, {\it Central charges in the canonical realization of asymptotic symmetries: an example from three-dimensional gravity}, Commun. Math. Phys. {\bf 104} (1986), 207--226.

\bibitem[CFP96]{CFP96} J. W. Cannon, W. J. Floyd, W. R. Parry, {\it Introductory notes on Richard Thompson's groups}, Enseign. Math. {\bf 42} (1996), 215--256.

\bibitem[CF97]{CF} L. O. Chekhov, V. V. Fock, {\it A quantum Teichm\"uller space}, Theor. Math. Phys. {\bf 120} (1999), 511-528.

\bibitem[F95]{Faddeev} L. D. Faddeev, {\it Discrete Heisenberg-Weyl group and modular group}, Lett. Math. Phys. {\bf 34} (1995), 249--254.

\bibitem[FK94]{FK94} L. D. Faddeev, R. M. Kashaev, {\it Quantum dilogarithm}, Mod. Phys. Lett. A {\bf 9} (1994), 427--434.

\bibitem[FKV01]{FKV} L. D. Faddeev, R. M. Kashaev, A. Y. Volkov, {\it Strongly coupled quantum discrete Liouville theory, I: Algebraic approach and duality}, Commun. Math. Phys. {\bf 219} (2001), 199-219.

\bibitem[F97]{F} V. V. Fock, {\it Dual Teichm\"uller spaces}, arXiv:dg-ga/9702018

\bibitem[FG06]{FG06} V. V. Fock, A. B. Goncharov, {\it Moduli spaces of local systems and higher Teichm\"uller theory}, Publ. Math. Inst. Hautes \'Etudes Sci. {\bf 103} (2006), 1--211.

\bibitem[FG09]{FG09} V. V. Fock, A. B. Goncharov, {\it The quantum dilogarithm and representations of the quantum cluster varieties}, Invent. Math. {\bf 175}(2) (2009), 223-286.

\bibitem[FST08]{FST} S. Fomin, M. Shapiro and D. Thurston, {\it Cluster algebras and triangulated surfaces. Part I: cluster complexes}, Acta Math. {\bf 201}(1) (2008), 83--146. \quad arXiv:math/0608367

\bibitem[FZ02]{FZ} S. Fomin, A. Zelevinsky, {\it Cluster algebras I: Foundations}, J. Amer. Math. Soc. {\bf 15} (2002), 497--529.

\bibitem[FK12]{FK} I. B. Frenkel, H. Kim, {\it Quantum Teichm\"uller space from the quantum plane}, Duke Math. J. {\bf 161}(2) (2012), 305--366.

\bibitem[FLM88]{FLM} I. B. Frenkel, J. Lepowsky, A. Meurman, {\it Vertex operator algebras and the Monster}, Academic Press, Boston, 1988.


\bibitem[FP22]{FP} I. B. Frenkel, R. C. Penner, {\it Sketch of a program for automorphic functions from universal Teichm\"uller theory to capture Monstrous Moonshine}, Commun. Math. Phys. {\bf 389} (2022), 1525--1567.

\bibitem[F95]{Funar} L. Funar, {\it $2+1$-D topological quantum field theory and $2$-D conformal field theory}, Commun. Math. Phys. {\bf 171} (1995), 405--458.


\bibitem[FS10]{FS} L. Funar, V. Sergiescu, {\it Central extensions of the Ptolemy-Thompson group and quantized Teichm\"uller theory}, J. Topol. {\bf 3} (2010), 29--62.

\bibitem[GTZ15]{GTZ} S. Garoufalidis, D. P. Thurston, C. K. Zickert, {\it The complex volume of ${\rm SL}(n,\mathbb{C})$-representations of $3$-manifold}, Duke Math. J. {\bf 164}(11) (2015), 2099--2160.

\bibitem[G90]{Goldman} W. Goldman, {\it Convex real projective structures on compact surfaces}, J. Differ. Geom. {\bf 31} (1990), 791--845.

\bibitem[G08]{G08} A. B. Goncharov, ``Pentagon relation for the quantum dilogarithm and quantized $\mathcal{M}^{\rm cyc}_{0,5}$" in {\it Geometry and Dynamics of Groups and Spaces}, Progr. Math. {\bf 265}, Birkh\"auser, Basel, 2008. pp. 415--428.

\bibitem[H13]{Hall} B. C. Hall, {\it Quantum theory for mathematicians}, Grad. Texts in Math., vol. 267, Springer, New York, 2013.



\bibitem[K98]{Kash98} R. M. Kashaev, {\it Quantization of Teichm\"uller spaces and the quantum dilogarithm}, Lett. Math. Phys. {\bf 43} (1998), 105--115.

\bibitem[K16a]{K16a} H. Kim, {\it The dilogarithmic central extension of the Ptolemy-Thompson group via the Kashaev quantization}, Adv. Math. {\bf 293} (2016), 529--588.

\bibitem[K16b]{K16b} H. Kim, {\it Ratio coordinates for higher Teichm\"uller spaces}, Math. Z. {\bf 283} (2016), 469--513.


\bibitem[K21a]{K21a} H. Kim, {\it Phase constants in the Fock-Goncharov quantum cluster varieties}, Anal. Math. Phys. {\bf 11} (2021), 2.

\bibitem[K21b]{K21b} H. Kim, {\it Irreducible self-adjoint representations of quantum Teichm\"uller space and the phase constants}, J. Geom. Phys. {\bf 162} (2021), 104103.

\bibitem[K23]{KimQC} H. Kim, {\it Three-dimensional quantum gravity from the quantum pseudo-K\"ahler plane}, Commun. Math. Phys. {\bf 402} (2023), 2715--2763. 

\bibitem[K24]{Kim24} H. Kim, {\it Quantized geodesic lengths for Teichmüller spaces: algebraic aspects}, arXiv:2405.14727

\bibitem[KS24]{KS} H. Kim, C. Scarinci, {\it A quantization of moduli spaces of $3$-dimensional gravity}, Commun. Math. Phys. {\bf 405} (2024), article number 144.

\bibitem[L-F09]{LF} D. Labardini-Fragoso, {\it Quivers with potentials associated to triangulated surfaces, part II: arc representations}, arXiv:math.RT/0909.4100

\bibitem[L77]{Lion} G. Lion, ``Integrales d'Entrelacement sur des groupes de Lie nilpotents et indices de Maslov" in J. Carmona, M. Vergne (Eds.), {\it Non-commutative harmonic analysis}, Lect. Notes in Math., vol. 587, Springer, Berlin, Heidelberg, 1977, pp. 160--176.

\bibitem[MS16]{MS} C. Meusburger, C. Scarinci, {\it Generalized shear coordinates on the moduli spaces of three-dimensional spacetimes}, J. Differ. Geom. {\bf 103} (2016), 425--474.

\bibitem[NN04]{NN} T. Nakanishi, M. Näätänen, {\it Complexification of lambda length as parameter for ${\rm SL}(2,\mathbb{C})$ representation space of punctured surface groups}, J. London Math. Soc. (2) {\bf 70} (2004), 383--404.


\bibitem[P87]{P87} R. C. Penner, {\it The decorated Teichm\"uller space of punctured surfaces}, Commun. Math. Phys. {\bf 113} (1987), 299-339.

\bibitem[P93]{Penner93} R. C. Penner, {\it Universal constructions in Teichm\"uller theory}, Adv. Math. {\bf 98} (1993), 143--215.

\bibitem[SK13]{SK} C. Scarinci, K. Krasnov, {\it The universal phase space of ${\rm AdS}_3$ gravity}, Commun. Math. Phys. {\bf 322} (2013), 167--205.

\bibitem[S12]{Sch} K. Schm\"udgen, {\it Unbounded self-adjoint operators on Hilbert space}, Grad. Texts in Math., vol. 265, Springer, Dordrecht, 2012.

\bibitem[S92]{Schpaper} K. Schm\"udgen, {\it Operator representations of $\mathbb{R}^2_q$}, Publ. Res. Inst. Math. Sci. {\bf 28} (1992), 1029--1061.

\bibitem[S63]{Segal} I. E. Segal, {\it Transforms for operators and symplectic automorphisms over a locally compact abelian group}, Math. Scand. {\bf 12}--{\bf 13} (1963), 31--43.

\bibitem[S62]{Shale} D. Shale, {\it Linear symmetries of free boson fields}, Trans. Amer. Math. Soc. {\bf 103} (1962), 149--167.

\bibitem[T07]{Teschner} J. Teschner, {\it An analog of a modular functor from quantized Teichm\"uller theory}, Handbook of Teichm\"uller theory {\bf Vol. I}, IRMA Lect. Math. Theor. Phys., {\bf 11} (2007), 685--760. 

\bibitem[T80]{Thurston} W.P. Thurston, ``The Geometry and Topology of three-manifolds" in: Lecture Notes, Princeton University, 1980, available at http://library.msri.org/books/gt3m

\bibitem[vN31]{vN} J. von Neumann, {\it Die Eindeutigkeit der Schr\"odingerschen operatoren}, Math. Ann. {\bf 104} (1931), 570--578.

\bibitem[W64]{Weil} A. Weil, {\it Sur certaines groupes d'op\'erateurs unitaires}, Acta Math. {\bf 11} (1964), 143--211.

\bibitem[W07]{Witten} E. Witten, {\it Three-dimensional gravity reconsidered}, arXiv:0706.3359.

\bibitem[W00]{W00} S. L. Woronowicz, {\it Quantum exponential function}, Rev. Math. Phys. {\bf 12}(6) (2000), 873--920.

\bibitem[Z07]{Zagier} D. Zagier, ``The dilogarithm function" in P. Cartier, P. Moussa, B. Julia, P. Vanhove (eds), {\it Frontiers in number theory, physics and geometry II}, Springer, Berlin, Heidelberg, 2007, pp.3--65.



\end{thebibliography}

\end{document}